%% file: main_arxiv.tex
\documentclass{article}
\usepackage[a4paper, margin=3cm]{geometry}
\usepackage{graphicx} 
\usepackage[linktocpage,colorlinks,linkcolor=blue,anchorcolor=blue,citecolor=blue,urlcolor=blue,pagebackref]{hyperref}       
\usepackage{url}           
\usepackage{booktabs}       
\usepackage{amsfonts}       
\usepackage{nicefrac}      
\usepackage{microtype}      
\usepackage{xcolor}   

\usepackage{amsmath}
\usepackage{amssymb}
\usepackage{multirow}
\usepackage{makecell}
\usepackage{booktabs}
\usepackage{array}
\usepackage{enumitem}

\usepackage{arydshln}
\usepackage{multicol}
\usepackage{algorithm}
\usepackage{algorithmic}
\usepackage{amsthm}

\usepackage{wrapfig}

\usepackage[capitalize,noabbrev]{cleveref}
\usepackage{breakcites}

\newtheorem{theorem}{Theorem}[section]

\newtheorem{remark}{Remark}
\newtheorem{definition}{Definition}[section]
\newtheorem{lemma}{Lemma}

\newtheorem{assumption}{Assumption}
\newenvironment{keywords}
{\par\smallskip\noindent\textbf{Keywords:}\ }
{\par\smallskip}
\input{macro}

\usepackage{natbib}

\usepackage[most]{tcolorbox}

\providecommand{\email}[1]{\href{mailto:#1}{\texttt{#1}}}
\title{Newton-Based Methods for Bilevel Optimization Beyond Lower-Level Strong Convexity}
\author{
	Youran Dong\thanks{School of Mathematics, Nanjing University, Nanjing, P. R. China. (\email{yrdong@smail.nju.edu.cn}).} \quad
	Junfeng Yang\thanks{School of Mathematics, Nanjing University, Nanjing, P. R. China. (\email{jfyang@nju.edu.cn}).} \quad
	Wei Yao\thanks{School of Artificial Intelligence, Southern University of Science and Technology, Shenzhen, P. R. China. (\email{yaow@sustech.edu.cn}).} \quad
	Jin Zhang\,\thanks{Department of Mathematics and National Center for Applied Mathematics Shenzhen, Southern University of Science and Technology, Shenzhen, P. R. China. (\email{zhangj9@sustech.edu.cn}).}
}

\date{}

\begin{document}

\maketitle
\begin{abstract}
	Bilevel optimization provides a framework for machine learning tasks such as hyperparameter optimization, data reweighting, and meta-learning. Existing Newton-based bilevel methods typically rely on the classical hypergradient formula, which is generally unavailable when lower-level solutions are nonunique and the lower-level Hessian is singular. We develop two Newton-based methods beyond lower-level strong convexity by exploiting curvature shared by the lower-level and multiplier updates. Newton-based GALET handles potentially singular lower-level curvature under a positive spectral-gap condition, exploiting its range structure to avoid the squared-Hessian multiplier system used in GALET. Newton-based BAMM uses a strongly convex aggregation model when the lower-level objective is convex and the upper-level objective is strongly convex in the lower-level variable. We establish convergence guarantees and improved complexity bounds relative to GALET and sl-BAMM under the corresponding assumptions. We further introduce GNBO, a generalized Newton-based framework that captures the shared-curvature structure of the proposed methods and several existing Newton and quasi-Newton bilevel methods. Experiments on synthetic problems, data hyper-cleaning, and GPT-2 data reweighting evaluate the proposed methods across different bilevel learning tasks.
\end{abstract}


\section{Introduction}

Bilevel optimization models hierarchical learning and decision-making problems in which one optimization problem depends on the solution of another. It arises in hyperparameter optimization \citep{pedregosa2016hyperparameter,franceschi2018bilevel}, data cleaning and reweighting \citep{ren2018learning,pan2025scalebio}, meta-learning \citep{franceschi2018bilevel,rajeswaran2019meta,ji2020convergence}, and adversarial learning \citep{bishop2020optimal,wang2021fast,zhang2022revisiting}. We consider the optimistic bilevel problem
\begin{equation}
\min_{x\in\mathbb{R}^m,\;y\in Y^*(x)} f(x,y),
\qquad
Y^*(x):=\arg\min_{y\in\mathbb{R}^n} g(x,y),
\label{eq:bo}
\end{equation}
where $f$ and $g$ are the upper-level (UL) and lower-level (LL) objectives, respectively, and the LL problem may admit multiple solutions.

One of the main challenges in solving \eqref{eq:bo} stems from its nested optimization structure. When $g(x,\cdot)$ is strongly convex, the LL problem admits a unique solution $y^*(x)$, and the bilevel problem reduces to minimizing the hyper-objective $\Phi(x):=f(x,y^*(x))$. Under standard smoothness assumptions, the implicit function theorem yields the hypergradient
\begin{equation}
\nabla\Phi(x)
=
\nabla_1 f\bigl(x,y^*(x)\bigr)
-
\nabla_{12}^2 g\bigl(x,y^*(x)\bigr)
\bigl[\nabla_{22}^2 g\bigl(x,y^*(x)\bigr)\bigr]^{-1}
\nabla_2 f\bigl(x,y^*(x)\bigr).
\label{eq:hypergradient}
\end{equation} 

A large body of bilevel optimization methods approximates the hypergradient in \eqref{eq:hypergradient} and applies a gradient-type update to the UL variable. A standard approach is approximate implicit differentiation (AID) \citep{pedregosa2016hyperparameter,ghadimi2018approximation,grazzi2020iteration,ji2022will,chen2021closing,hong2023two,dagreou2022framework,chu2024spaba,ramzi2022shine,fang2025qnbo,dong2025efficient}. AID involves two key approximation tasks: estimating the LL solution $y^*(x)$ and approximating the Hessian-inverse vector product $\bigl[\nabla_{22}^2 g\bigl(x,y^*(x)\bigr)\bigr]^{-1}
\nabla_2 f\bigl(x,y^*(x)\bigr)$. Importantly, the Hessian $\nabla_{22}^2 g$ in the second task is the curvature of the same LL objective whose minimizer is estimated in the first, revealing a shared-curvature structure between the two tasks. 

Motivated by this observation, several recent Newton-based bilevel methods exploit LL curvature to couple these two approximation tasks in different ways. SHINE \citep{ramzi2022shine} applies a quasi-Newton method to the LL problem and reuses the resulting quasi-Newton approximation to estimate the Hessian inverse in the hypergradient. qNBO \citep{fang2025qnbo} instead uses separate quasi-Newton procedures for the LL solution and the Hessian-inverse vector product, with nonasymptotic guarantees based on properties of BFGS updates. More recently, NBO \citep{dong2025efficient} uses two inexact Newton directions involving the same LL Hessian to approximate the LL solution and the Hessian-inverse vector product, leading to improved nonasymptotic complexity bounds over existing AID methods. Together, these methods highlight the benefits of exploiting LL curvature in both LL optimization and hypergradient approximation. 

The Newton-based methods discussed above rely on LL strong convexity and the resulting hypergradient formula \eqref{eq:hypergradient}. When the LL problem admits multiple solutions, however, the solution mapping becomes set-valued and the implicit differentiation underlying \eqref{eq:hypergradient} is generally unavailable. An alternative is to characterize LL optimality through first-order stationarity and consider
\begin{equation}
\min_{x\in\mathbb{R}^m,\;y\in\mathbb{R}^n} f(x,y)
\qquad
\mathrm{s.t.}
\qquad
\nabla_2 g(x,y)=0.
\label{eq:boeq}
\end{equation}
When $g(x,\cdot)$ is convex or satisfies the Polyak--{\L}ojasiewicz (PL) condition, its stationary points are global minimizers, making \eqref{eq:boeq} equivalent to the original bilevel problem \eqref{eq:bo}. This reformulation underlies several methods for bilevel optimization beyond LL strong convexity
\citep{mehra2021penalty,shen2025penalty,xiao2023generalized,liu2023averaged,hallak2026augmented}. 

The constrained viewpoint has led to several methods beyond LL strong convexity. \cite{xiao2023generalized} derive KKT-type necessary conditions for bilevel problems with PL LL objectives and develop GALET, whose convergence analysis assumes a uniform positive gap in the nonzero singular values of the LL Hessian. \cite{liu2023averaged} introduce an iteration-dependent aggregation model and develop the single-loop method sl-BAMM to approximate the resulting KKT-type residuals. Fully first-order alternatives include BOME \citep{liu2022bome}, which uses a value-function-based constrained reformulation and a modified dynamic barrier gradient method, and HJFBiO \citep{huang2024optimal}, which combines finite-difference estimators with projection-based updates for Hessian/Jacobian-free nonconvex-PL bilevel optimization under an additional nonsingularity condition on the LL Hessian at LL solutions. 


\subsection{Motivation and Contributions}

The preceding discussion motivates extending Newton-type methods beyond LL strong convexity. When LL solutions are nonunique, the LL Hessian may be singular, making the classical hypergradient formula generally unavailable. This raises a natural question:

\begin{center}
	\textit{Can Newton-type methods effectively address bilevel optimization \\ with possibly nonunique lower-level solutions?}
\end{center}

Our main contributions are summarized as follows.
\begin{enumerate}[label=(\roman*)]
	\item We develop \emph{Newton-based GALET} for bilevel problems with possibly nonunique LL solutions under a positive spectral-gap condition, as also assumed in GALET \citep{xiao2023generalized}. Starting from the KKT residuals of \eqref{eq:boeq}, we establish a range-inclusion property of the LL Hessian and exploit the curvature shared by the LL and multiplier residuals to construct Newton-type updates. This structure allows the multiplier update to use the LL Hessian directly, rather than the squared-Hessian system employed by GALET. 
	
	\item We develop \emph{Newton-based BAMM} for problems with a convex LL objective and an UL objective strongly convex in the LL variable. Building on the aggregation model of sl-BAMM \citep{liu2023averaged}, we construct Newton-type LL and multiplier updates using the same positive-definite Hessian.
	
	\item We establish convergence and complexity guarantees for both methods under assumptions comparable to those of GALET \citep{xiao2023generalized} and sl-BAMM \citep{liu2023averaged}, respectively. The resulting bounds improve the corresponding complexity results, with gradient/Jacobian-vector-product (JVP) and Hessian-vector-product costs analyzed separately. We further introduce GNBO, a generalized Newton-based framework that captures the shared-curvature structure of the proposed methods and several existing Newton and quasi-Newton approaches \citep{ramzi2022shine,fang2025qnbo,dong2025efficient}.
	
	\item We evaluate the proposed methods on two synthetic examples, data hyper-cleaning with linear and neural models, and GPT-2 data reweighting, examining their empirical performance across different bilevel learning tasks.
\end{enumerate}


\subsection{Notation}

For a differentiable function $g$, we denote its partial gradients with respect to $x$ and $y$ by $\nabla_1g(x,y)$ and $\nabla_2g(x,y)$, respectively. For twice differentiable $g$, $\nabla_{12}^2g(x,y)\in\mathbb{R}^{m\times n}$ denotes the Jacobian of $\nabla_1g(x,y)$ with respect to $y$, and $\nabla_{22}^2g(x,y)\in\mathbb{R}^{n\times n}$ denotes the Hessian with respect to $y$. Unless otherwise specified, $\|\cdot\|$ denotes the Euclidean norm for vectors and the Frobenius norm for matrices, while $\|\cdot\|_{\mathrm{op}}$ denotes the spectral norm. 

\subsection{Organization}

The remainder of the paper is organized as follows. Section~\ref{methods} develops Newton-based GALET and Newton-based BAMM. Sections~\ref{analysis} and \ref{experiments} present the convergence analysis and numerical experiments, respectively. Section~\ref{sec:discussion} introduces GNBO as a unifying framework, and Section~\ref{conclude} concludes the paper. Detailed proofs and experimental details are provided in the appendices. 

\section{Newton-Based Bilevel Optimization Beyond Strong Convexity}
\label{methods}

\subsection{From Hypergradients to KKT Residuals}

We first revisit the shared-curvature structure in the strongly convex LL setting. When $g(x,\cdot)$ is strongly convex, the LL solution $y^*(x)$ is unique, and the hypergradient error induced by an inexact LL solution $y$ is
$O(\|y-y^*(x)\|)$ under standard smoothness assumptions
\citep{ghadimi2018approximation}. Locally, a Newton step
$y^+=y-[\nabla_{22}^2g(x,y)]^{-1}\nabla_2g(x,y)$
reduces the LL solution error quadratically:
$\|y^+-y^*(x)\|=O(\|y-y^*(x)\|^2)$.
The Hessian-inverse vector product in \eqref{eq:hypergradient} involves the same LL Hessian. Exploiting this shared-curvature structure, NBO \citep{dong2025efficient} uses the following two-step procedure:
\begin{description}
	\item[(i) \textbf{Newton step.}]
	Approximately solve
	$\nabla_{22}^2g(x,y)(v,u)
	=(\nabla_2g(x,y),\nabla_2f(x,y))$.
	
	\item[(ii) \textbf{Hypergradient step.}]
	Set
	$d_x:=\nabla_1f(x,y-v)-\nabla_{12}^2g(x,y-v)u$
	and update the UL variable.
\end{description}
Here, $-v$ is an inexact Newton direction for the LL variable, while $u$ approximates the Hessian-inverse vector product. Both involve the same LL Hessian. 

When the LL solution may be nonunique, the classical hypergradient formula
\eqref{eq:hypergradient} is generally unavailable. We instead consider
\eqref{eq:boeq}, with Lagrangian
$\mathcal{L}(x,y,u)=f(x,y)-u^\top\nabla_2g(x,y)$
and KKT residuals
\begin{align}
r_y(x,y)&=\nabla_2g(x,y), \label{eq:kkt1}\\
r_u(x,y,u)&=\nabla_2f(x,y)-\nabla_{22}^2g(x,y)u, \label{eq:kkt2}\\
r_x(x,y,u)&=\nabla_1f(x,y)-\nabla_{12}^2g(x,y)u. \label{eq:kkt3}
\end{align}
Under the PL condition and suitable regularity assumptions,
$r_y=r_u=r_x=0$ is a necessary optimality condition
\citep{xiao2023generalized}. Under strong convexity, these conditions recover the LL solution, the Hessian-inverse vector product, and hypergradient stationarity, respectively \citep{liu2023averaged}.

The shared-curvature structure extends naturally to these residuals. For fixed $x$, Taylor expansion gives
\begin{equation}
r_y(x,y-v)
=
r_y(x,y)-\nabla_{22}^2g(x,y)v+O(\|v\|^2).
\label{eq:kkty1}
\end{equation}
Hence, if
$\nabla_{22}^2g(x,y)v=\nabla_2g(x,y)$
and $\|v\|=O(\|\nabla_2g(x,y)\|)$, then
$\|\nabla_2g(x,y-v)\|
=O(\|\nabla_2g(x,y)\|^2)$.
For the multiplier residual, $r_u(x,y,u)=0$ corresponds to
$\nabla_{22}^2g(x,y)u=\nabla_2f(x,y)$.
Thus, the LL and multiplier residuals give rise to two systems involving the same LL Hessian. 

When the LL solution is nonunique, this Hessian may be singular, so the two systems need not have unique solutions. Rather than relying on a Hessian inverse, we therefore work directly with their residuals:
\begin{description}
	\item[(i) \textbf{Residual step.}]
	Given $(x,y)$, construct $v$ and $u$ by reducing
	$r_1:=\nabla_2g(x,y)-\nabla_{22}^2g(x,y)v$
	and
	$r_2:=\nabla_2f(x,y)-\nabla_{22}^2g(x,y)u$.
	
	\item[(ii) \textbf{UL step.}]
	Construct $d_x$ from $r_x(x,y-v,u)$ and use it to update $x$.
\end{description}
Here, $r_1$ is the linearized residual associated with $r_y(x,y-v)$, so $-v$ is a Newton-type direction, while $r_2=r_u(x,y,u)$ is the multiplier residual. This residual viewpoint preserves the shared-curvature structure without requiring an invertible LL Hessian and allows the two systems to be treated according to the problem structure. 

We next instantiate this principle in two settings, leading to Newton-based GALET and Newton-based BAMM.


\subsection{Newton-Based GALET}

We first consider the positive spectral-gap setting.

\begin{assumption}\label{asp:singularvalue}
	There exists $\mu>0$ such that every nonzero singular value of
	$\nabla_{22}^2g(x,y)$ is uniformly bounded below by $\mu$, i.e.,
	\begin{equation*}
	\inf_{x,y}\sigma_{\min}^{+}\bigl(\nabla_{22}^2g(x,y)\bigr)\geq\mu,
	\end{equation*}
	where $\sigma_{\min}^{+}$ denotes the smallest positive singular value.
\end{assumption}

Assumption~\ref{asp:singularvalue}, also used in GALET
\citep{xiao2023generalized}, allows the LL Hessian to be singular and is therefore compatible with nonunique LL solutions. For example,
$g(x,y)=\frac{1}{2}\|A(x-y)-b\|^2$
satisfies this assumption for a fixed matrix
$A\in\mathbb{R}^{n\times n}$ with $\rank A=r<n$.
It also implies
$\|[\nabla_{22}^2g(x,y)]^\dagger\|_{\operatorname{op}}\leq1/\mu$.
The following lemma characterizes the resulting LL structure. Its proof is given in Appendix~\ref{apd:proof}. 

\begin{lemma}\label{lem:range}
	Suppose Assumption~\ref{asp:singularvalue} holds,
	$g^*(x):=\min_y g(x,y)>-\infty$, and
	$\nabla_{22}^2g(x,\cdot)$ is $L_{g,2}$-Lipschitz. Then $g(x,\cdot)$ is convex,
	$\nabla_2g(x,y)\in\range\He$, and
	$\|\nabla_2g(x,y)\|\geq
	\mu\operatorname{dist}(y,Y^*(x))$. 
	If, in addition, $\nabla_2g(\cdot,y)$ is $L_{g,1}$-Lipschitz, then
	$\range{\nabla_{21}^2g(x,y)}\subseteq\range\He$.
\end{lemma}

\begin{remark}
	By convexity and the error bound in Lemma~\ref{lem:range},
	\begin{equation*}
	g(x,y)-g^*(x)
	\leq \|\nabla_2g(x,y)\|\operatorname{dist}(y,Y^*(x))
	\leq \mu^{-1}\|\nabla_2g(x,y)\|^2.
	\end{equation*}
	Thus, Assumption~\ref{asp:singularvalue}, together with the conditions of Lemma~\ref{lem:range}, implies the PL condition with parameter $\mu/2$, which is assumed separately in GALET \citep{xiao2023generalized}.
\end{remark} 

\paragraph{Algorithm Design}
To reduce the LL residual $r_y$, we approximately solve
\begin{equation}
\Hk v=\nabla_2g(x^k,y^k)
\label{eq:galet-v-system}
\end{equation}
for $v^k$ and set $y^{k+1}=y^k-v^k$. Although $\Hk$ may be singular, Lemma~\ref{lem:range} ensures
$\nabla_2g(x^k,y^k)\in\range\Hk$,
so \eqref{eq:galet-v-system} is consistent. 
We approximately solve the linear systems using Algorithm~\ref{alg:PLGD},
where $\mathcal{H}$ denotes a linear operator. Throughout, Hessian-vector products
are computed via automatic differentiation
\citep{pearlmutter1994fast,dagreou2024hvp},
without explicitly forming or storing the Hessian matrices. 

\begin{algorithm}[ht]
	\caption{$\operatorname{GD}(\mathcal{H},q;\gamma,T)$ for $Hp=q$, with $\mathcal{H}(p)=Hp$}
	\label{alg:PLGD}
	\begin{algorithmic}[1]
		\STATE Initialize $p^0=0$.
		\FOR{$t=0,\ldots,T-1$}
		\STATE $p^{t+1}=p^t-\gamma(\mathcal{H}(p^t)-q)$.
		\ENDFOR
		\STATE \textbf{Return} $p^T$.
	\end{algorithmic}
\end{algorithm} 

For the multiplier update, let
$q_k:=\nabla_2f(x^k,y^k)$, so
$r_2=q_k-\Hk u$.
Since $q_k$ need not lie in $\range\Hk$, $r_2=0$ may be infeasible.
Let $P_k:=\Hk[\Hk]^\dagger$ be the orthogonal projector onto
$\range\Hk$. Then
\begin{equation*}
r_2
=
(P_kq_k-\Hk u)+(I-P_k)q_k,
\end{equation*}
where the two components are orthogonal and the second is independent of $u$.
We therefore reduce the controllable component through the consistent system
\begin{equation}
\Hk u=P_kq_k.
\label{eq:galet-u-system}
\end{equation}

The projection need not be computed explicitly. We justify this for both the multiplier residual and the UL direction.

First, applying Algorithm~\ref{alg:PLGD} directly to $(\Hk,q_k)$ gives
\[
p^{t+1}=p^t-\eta_k(\Hk p^t-q_k).
\] 
Since $\Hk$ is symmetric,
$P_k\Hk=\Hk P_k=\Hk$, and hence
\begin{equation*}
P_kp^{t+1}
=
P_kp^t-\eta_k(\Hk P_kp^t-P_kq_k).
\end{equation*}
Thus, $\{P_kp^t\}$ is exactly the GD sequence for the consistent system
\eqref{eq:galet-u-system}.

Second, Lemma~\ref{lem:range} implies
$\nabla_{12}^2g(x^k,y^k)(I-P_k)=0$.
Consequently,
\[
\nabla_{12}^2g(x^k,y^k)P_kp^t
=\nabla_{12}^2g(x^k,y^k)p^t,\qquad
\Hk P_kp^t=\Hk p^t.
\]
Hence, the component of $p^t$ in $\ker\Hk$ affects neither the UL direction nor the Hessian-weighted multiplier residual. We can therefore set
$u^{k+1}=\operatorname{GD}(\Hk,q_k;\eta_k,T_2)$
without explicitly forming $P_k$. 

We then set
\begin{align*}
d_x^k
=
\nabla_1f(x^k,y^k)
-\nabla_{12}^2g(x^k,y^k)u^{k+1},
\end{align*}
and update $x^{k+1}=x^k-\alpha_kd_x^k$.
Algorithm~\ref{alg:ngalet} summarizes Newton-based GALET. Unlike GALET \citep{xiao2023generalized}, which solves a squared-Hessian system, Newton-based GALET solves the LL Hessian system directly, requiring only one Hessian-vector product per inner GD step instead of two.

\begin{algorithm}[h]
	\caption{Newton-Based GALET}
	\label{alg:ngalet}
	\begin{algorithmic}[1]
		\STATE \textbf{Input:} $x^0,y^0$; $K,T_1,T_2$; step sizes $\{\alpha_k,\gamma_k,\eta_k\}$.
		\FOR{$k=0,\ldots,K-1$}
		\STATE Define $\mathcal{H}_k(p):=\nabla_{22}^2g(x^k,y^k)p$.
		\STATE $v^k=\operatorname{GD}(\mathcal{H}_k,\nabla_2g(x^k,y^k);\gamma_k,T_1)$.
		\STATE $u^{k+1}=\operatorname{GD}(\mathcal{H}_k,\nabla_2f(x^k,y^k);\eta_k,T_2)$.
		\STATE Update $y^{k+1}=y^k-v^k$.
		\STATE Set $d_x^k=\nabla_1f(x^k,y^k)-\nabla_{12}^2g(x^k,y^k)u^{k+1}$.
		\STATE Update $x^{k+1}=x^k-\alpha_kd_x^k$.
		\ENDFOR
	\end{algorithmic}
\end{algorithm}


\subsection{Newton-Based BAMM}

We next consider problems with a convex LL objective and an UL objective strongly convex in the LL variable.

\begin{assumption}\label{asp:fsc}
	The lower-level function $g(x,y)$ is convex in $y$, and the upper-level function $f(x,y)$ is $\mu_f$-strongly convex in $y$.
\end{assumption}

sl-BAMM \citep{liu2023averaged} uses the aggregation function
$\psi_{\sigma_k}(x,y):=\sigma_k f(x,y)+(1-\sigma_k)g(x,y)$,
where $0<\sigma_k<1$ decreases to zero. By Assumption~\ref{asp:fsc},
$\psi_{\sigma_k}(x,\cdot)$ is $(\sigma_k\mu_f)$-strongly convex, with unique minimizer
$y_{\sigma_k}^*(x):=\arg\min_y\psi_{\sigma_k}(x,y)$.
The corresponding residuals are
\begin{align*}
r_y(x,y)&=\nabla_2\psi_{\sigma_k}(x,y),\\
r_u(x,y,u)&=\nabla_2f(x,y)-\nabla_{22}^2\psi_{\sigma_k}(x,y)u,\\
r_x(x,y,u)&=\nabla_1f(x,y)-\nabla_{12}^2\psi_{\sigma_k}(x,y)u.
\end{align*}
As $\sigma_k\to0$, these residuals approach those in
\eqref{eq:kkt1}--\eqref{eq:kkt3}. The LL and multiplier residuals share the aggregation Hessian
$\nabla_{22}^2\psi_{\sigma_k}$. 

\paragraph{Algorithm Design}
At iteration $k$, we approximately solve
\begin{equation}
\nabla_{22}^2\psi_{\sigma_k}(x^k,y^k)v
=
\nabla_2\psi_{\sigma_k}(x^k,y^k)
\label{eq:bamm-v-system}
\end{equation}
for $v^k$ and set $y^{k+1}=y^k-v^k$. For the multiplier, we use the warm-start update
$u^{k+1}=u^k-w^k$, where $w^k$ approximately solves
\begin{equation}
\nabla_{22}^2\psi_{\sigma_k}(x^k,y^k)w
=
\nabla_{22}^2\psi_{\sigma_k}(x^k,y^k)u^k
-\nabla_2f(x^k,y^k).
\label{eq:bamm-u-system}
\end{equation}
An exact solution of \eqref{eq:bamm-u-system} yields
$r_u(x^k,y^k,u^{k+1})=0$.

The two systems share the same aggregation Hessian. Unlike Newton-based GALET, no range-space projection is needed because
$\psi_{\sigma_k}(x^k,\cdot)$ is $(\sigma_k\mu_f)$-strongly convex, so its Hessian is positive definite for $\sigma_k>0$. As $\sigma_k\to0$, this modulus vanishes, and its decay must be controlled in the convergence analysis; see Theorem~\ref{thm:bamm}.

Following the UL residual, we set
$d_x^k=\nabla_1f(x^k,y^k)
-\nabla_{12}^2\psi_{\sigma_k}(x^k,y^k)u^k$
and update $x^{k+1}=x^k-\alpha_kd_x^k$.
Algorithm~\ref{alg:nbamm} summarizes Newton-based BAMM.

\begin{algorithm}[h] 
	\caption{Newton-Based BAMM}
	\label{alg:nbamm}
	\begin{algorithmic}[1]
		\STATE \textbf{Input:} $x^0,y^0,u^0$; $K,\{T_k\}_{k=0}^{K-1}$; step sizes $\{\alpha_k,\gamma_k\}$; aggregation parameters $\{\sigma_k\}$.
		\FOR{$k=0,\ldots,K-1$}
		\STATE Define $\mathcal{H}_k(p):=\nabla_{22}^2\psi_{\sigma_k}(x^k,y^k)p$.
		\STATE $v^k=\operatorname{GD}(\mathcal{H}_k,\nabla_2\psi_{\sigma_k}(x^k,y^k);\gamma_k,T_k)$.
		\STATE $w^k=\operatorname{GD}(\mathcal{H}_k,\mathcal{H}_k(u^k)-\nabla_2f(x^k,y^k);\gamma_k,T_k)$.
		\STATE Update $y^{k+1}=y^k-v^k$ and $u^{k+1}=u^k-w^k$.
		\STATE Set $d_x^k=\nabla_1f(x^k,y^k)-\nabla_{12}^2\psi_{\sigma_k}(x^k,y^k)u^k$.
		\STATE Update $x^{k+1}=x^k-\alpha_kd_x^k$.
		\ENDFOR
	\end{algorithmic}
\end{algorithm} 


\section{Convergence Analysis}
\label{analysis}

We establish convergence and complexity guarantees for Newton-based GALET and Newton-based BAMM; proofs are given in Appendix~\ref{apd:proof}. Both analyses use a Lyapunov argument. For Newton-based GALET, we combine the Lagrangian-type potential
$\mathcal{P}_k=f(x^k,y^k)-u^*(x^k,y^k)^\top\nabla_2g(x^k,y^k)$,
where
$u^*(x^k,y^k)=[\Hk]^\dagger\nabla_2f(x^k,y^k)$,
with the LL residual $\|\nabla_2g(x^k,y^k)\|^2$. For Newton-based BAMM, we combine
$\mathcal{P}_k=f(x^k,y_{\sigma_k}^*(x^k))$
with the LL and multiplier tracking errors.

The key ingredient is the local contraction induced by the inexact Newton LL update: the LL residual for Newton-based GALET and the LL tracking error for Newton-based BAMM, up to the UL perturbation and, for BAMM, the variation of $\sigma_k$. Suitable initialization and parameter choices keep the iterates in the required local regions. 


\subsection{Convergence Analysis of Newton-Based GALET}

In addition to Assumption~\ref{asp:singularvalue}, we impose the following regularity conditions.

\begin{assumption}\label{asp:smooth1}
	\begin{enumerate}[label=(\alph*)]
		\item The UL function $f(x,y)$ has an $L_{f,1}$-Lipschitz continuous gradient.
		\item The LL function $g(x,y)$ has an $L_{g,1}$-Lipschitz continuous gradient and an $L_{g,2}$-Lipschitz continuous Hessian.
		\item There exists $C_{f,0}>0$ such that
		$\|\nabla_1f(x,y)\|\leq C_{f,0}$ for all $x,y$.
		\item There exists $L_{f,0}>0$ such that
		$\|\nabla_2f(x,y)\|\leq L_{f,0}$ for all $x,y$.
	\end{enumerate}
\end{assumption}

\begin{remark}\label{rem:asp}
	Except for part~(c), these conditions are also assumed in
	\cite{xiao2023generalized}. Part~(c) bounds the UL direction and helps keep the LL iterates in the local region required for the Newton-type error reduction. It holds, for example, in data hyper-cleaning when the validation objective is independent of the UL variable, so $\nabla_1f(x,y)=0$.
\end{remark} 

Following the necessary optimality conditions in
\cite[Theorem~1]{xiao2023generalized}, we adopt the stationarity measure used in GALET.

\begin{definition}\label{def:galet-stationarity}
	A point $(\bar{x},\bar{y})$ is $\epsilon$-stationary if there exists $\bar{u}$ such that
	$\mathcal{R}_x(\bar{x},\bar{y},\bar{u})\leq\epsilon$,
	$\mathcal{R}_u(\bar{x},\bar{y},\bar{u})\leq\epsilon$, and
	$\mathcal{R}_y(\bar{x},\bar{y})\leq\epsilon$, where
	\begin{align*}
	\mathcal{R}_x(x,y,u)
	&:=\|\nabla_1f(x,y)-\nabla_{12}^2g(x,y)u\|^2,\\
	\mathcal{R}_u(x,y,u)
	&:=\|\nabla_{22}^2g(x,y)
	(\nabla_2f(x,y)-\nabla_{22}^2g(x,y)u)\|^2,\\
	\mathcal{R}_y(x,y)
	&:=g(x,y)-\min_z g(x,z).
	\end{align*}
\end{definition}

Here, $\mathcal{R}_u$ is a Hessian-projected multiplier residual,
accommodating a possibly singular LL Hessian.
Let $L:=\max\{L_{f,1},L_{g,1},L_{g,2}\}$ and
$\kappa:=L/\mu$. Without loss of generality, we assume
$0<\mu\leq1$.
We next establish the complexity of Newton-based GALET. 

\begin{theorem}\label{thm:galet}
	Suppose Assumptions~\ref{asp:singularvalue} and \ref{asp:smooth1} hold, and
	$f^*:=\inf_{x,\,y\in Y^*(x)}f(x,y)>-\infty$.
	For Newton-based GALET, choose
	$T_1=\Theta(\kappa)$,
	$T_2=\Theta(\kappa\log(\kappa K))$,
	$\gamma_k,\eta_k=\Theta(1/L_{g,1})$ with
	$0<\gamma_k,\eta_k\leq1/L_{g,1}$,
	and $\alpha_k=\alpha=c_\alpha\kappa^{-4}$
	for a sufficiently small constant $c_\alpha>0$.
	The constants in $T_1$ and $T_2$ are chosen sufficiently large. 
	Suppose further that
	$\|\nabla_2g(x^0,y^0)\|
	\leq\min\{\mu^2/(2L_{g,2}),\mu^2/L\}$.
	Then, 
	\begin{align*}
	&\frac{1}{K}\sum_{k=0}^{K-1}\mathcal{R}_x(x^k,y^k,u^{k+1})
	=O(\kappa^4/K),\qquad
	\frac{1}{K}\sum_{k=0}^{K-1}\mathcal{R}_u(x^k,y^k,u^{k+1})
	=O(1/K),\\
	&\frac{1}{K}\sum_{k=0}^{K-1}\mathcal{R}_y(x^k,y^k)
	=O(1/(\kappa^2K)),
	\qquad \forall\, K\geq1. 
	\end{align*}
	Consequently, Newton-based GALET finds an $\epsilon$-stationary point in
	$O(\kappa^4/\epsilon)$ outer iterations, requiring
	$O(\kappa^4/\epsilon)$ gradient/JVP computations and
	$O\bigl(\kappa^5\epsilon^{-1}\log(\kappa/\epsilon)\bigr)$ Hessian-vector products.
\end{theorem}

\begin{remark}\label{initialization_discussion}
	The initialization condition in Theorem~\ref{thm:galet}
	can be satisfied by gradient descent on $g(x^0,\cdot)$,
	requiring $O(\kappa\log\kappa)$ LL gradient evaluations,
	with dependence on the initial optimality gap and other
	problem constants omitted. This is a one-time cost.
	Each subsequent outer iteration requires $O(1)$ gradient/JVP
	computations and $O(T_1+T_2)$ Hessian-vector products,
	yielding the stated oracle-complexity bounds.
\end{remark} 

Under comparable assumptions, Newton-based GALET improves the
gradient/JVP complexity of GALET from
$O(\kappa^6/\epsilon)$ to $O(\kappa^4/\epsilon)$
and the Hessian-vector-product complexity from
$O(\kappa^8\epsilon^{-1}\log(\kappa/\epsilon))$ to
$O(\kappa^5\epsilon^{-1}\log(\kappa/\epsilon))$.
Two features contribute to these improvements.
First, the inexact Newton LL update yields a local quadratic
reduction of the LL residual up to the linear-system error
(Lemma~\ref{lem:y}).
Second, the range inclusion
$\range{\nabla_{21}^2g(x,y)}
\subseteq\range{\nabla_{22}^2g(x,y)}$
allows the multiplier update to use the LL Hessian directly,
rather than the squared-Hessian system used in GALET. 

Compared with the GALET analysis \citep{xiao2023generalized}, two technical differences are worth noting: \textbf{(a)} We establish a global lower bound for the Lyapunov function using the LL error bound and the lower boundedness of the bilevel objective, enabling the telescoping argument underlying our complexity guarantees. \textbf{(b)} Our stationarity bounds are evaluated at $(x^k,y^k,u^{k+1})$ rather than $(x^k,y^k,u^k)$, as the cold-start multiplier update directly controls the Hessian-weighted multiplier residual at $u^{k+1}$ without requiring a uniform bound on $\{u^k\}$. 


\subsection{Convergence Analysis of Newton-Based BAMM}

In addition to Assumption~\ref{asp:fsc}, we impose the following regularity conditions.

\begin{assumption}\label{asp:smooth2}
	\begin{enumerate}[label=(\alph*)]
		\item The UL function $f(x,y)$ has an $L_{f,1}$-Lipschitz continuous gradient and an $L_{f,2}$-Lipschitz continuous Hessian, and is bounded below by $f_0$.
		\item The LL function $g(x,y)$ has an $L_{g,1}$-Lipschitz continuous gradient and an $L_{g,2}$-Lipschitz continuous Hessian.
		\item There exists $C_{f,0}>0$ such that
		$\|\nabla_1f(x,y)\|\leq C_{f,0}$ for all $x,y$.
	\end{enumerate}
\end{assumption}

Except for part~(c), these conditions are also assumed in
\cite{liu2023averaged}. As in the Newton-based GALET analysis, part~(c) is used to preserve the local region required for the Newton-type error reduction; see Remark~\ref{rem:asp}.

Under Assumptions~\ref{asp:fsc} and \ref{asp:smooth2},
$\psi_{\sigma_k}(x,\cdot)$ is $(\sigma_k\mu_f)$-strongly convex, with
$L_{\psi,1}$-Lipschitz gradient and $L_{\psi,2}$-Lipschitz Hessian, where
$L_{\psi,1}$ and $L_{\psi,2}$ are defined in \eqref{eq:constant-phi}.
Define
$y_\sigma^*(x):=\arg\min_y\psi_\sigma(x,y)$
and
$u_\sigma^*(x):=
[\nabla_{22}^2\psi_\sigma(x,y_\sigma^*(x))]^{-1}
\nabla_2f(x,y_\sigma^*(x))$.

Following \cite{liu2023averaged}, we call $(x,y,u)$ $\epsilon$-stationary if
$\operatorname{KKT}(x,y,u)\leq\epsilon$, where
\begin{equation*}
\operatorname{KKT}(x,y,u):=
\|\nabla_1f(x,y)-\nabla_{12}^2g(x,y)u\|^2+
\|\nabla_2f(x,y)-\nabla_{22}^2g(x,y)u\|^2+
\|\nabla_2g(x,y)\|^2.
\end{equation*}
As in sl-BAMM \citep{liu2023averaged}, we analyze three boundedness settings. Throughout this subsection, $O(\cdot)$ and $\Theta(\cdot)$ suppress dependence on
$\mu_f,L_{f,1},L_{f,2},L_{g,1}$, and $L_{g,2}$. 

\begin{theorem}\label{thm:bamm}
	Suppose Assumptions~\ref{asp:fsc} and \ref{asp:smooth2} hold, and let
	$0<\tau<1/3$. For Newton-based BAMM, choose
	$\sigma_k=\bar{\sigma}(k+k_0)^{-p}$ with
	$\bar{\sigma}>0$ and $0<\sigma_k<1$,
	$\gamma_k=\Theta(1/L_{\psi,1})$ with
	$0<\gamma_k\leq1/L_{\psi,1}$, and
	$T_k=\Theta(1/\sigma_k)$ sufficiently large so that
	$T_k\geq\log(1/4)/\log(1-\mu_f\sigma_k\gamma_k)$.
	Let $\bar{\alpha}>0$ be fixed, choose $k_0$ sufficiently large,
	and suppose
	$\|y^0-y_{\sigma_0}^*(x^0)\|
	\leq\mu_f\sigma_0/(2L_{\psi,2})$.
	Then the following results hold.
	
	\begin{itemize}
		\item \textbf{Case 1:}
		Suppose $\{x^k\}$, $\{y^k\}$, and $\{u^k\}$ are bounded.
		For $0<p<(1-3\tau)/6$, choose
		$\alpha_k=\bar{\alpha}(k+k_0)^{-4p-2\tau}$.
		Then
		\begin{equation*}
		\min_{0\leq k\leq K-1}
		\operatorname{KKT}(x^k,y^k,u^k)
		=
		O\left(K^{-(1-6p-3\tau)}+K^{-2p}\right).
		\end{equation*}
		Choosing $p=(1-3\tau)/8$ yields an $\epsilon$-stationary iterate using
		$O(\epsilon^{-4/(1-3\tau)})$ gradient/JVP computations and
		$O(\epsilon^{-(9-3\tau)/(2-6\tau)})$ Hessian-vector products.
		
		\item \textbf{Case 2:}
		Suppose
		$\{\nabla_2f(x^k,y_{\sigma_k}^*(x^k))\}$ and $\{u^k\}$ are bounded.
		For $0<p<(1-3\tau)/4$, choose
		$\alpha_k=\bar{\alpha}(k+k_0)^{-3p-2\tau}$.
		Then
		\begin{equation*}
		\min_{0\leq k\leq K-1}
		\operatorname{KKT}(x^k,y^k,u^k)
		=
		O\left(K^{-(1-4p-3\tau)}+K^{-2p}\right).
		\end{equation*}
		Choosing $p=(1-3\tau)/6$ yields an $\epsilon$-stationary iterate using
		$O(\epsilon^{-3/(1-3\tau)})$ gradient/JVP computations and
		$O(\epsilon^{-(7-3\tau)/(2-6\tau)})$ Hessian-vector products.
		
		\item \textbf{Case 3:}
		Suppose $\{u_{\sigma_k}^*(x^k)\}$ and $\{u^k\}$ are bounded.
		For $0<p<(1-3\tau)/2$, choose
		$\alpha_k=\bar{\alpha}(k+k_0)^{-2p-2\tau}$.
		Then
		\begin{equation*}
		\min_{0\leq k\leq K-1}
		\operatorname{KKT}(x^k,y^k,u^k)
		=
		O\left(K^{-(1-2p-3\tau)}+K^{-2p}\right).
		\end{equation*}
		Choosing $p=(1-3\tau)/4$ yields an $\epsilon$-stationary iterate using
		$O(\epsilon^{-2/(1-3\tau)})$ gradient/JVP computations and
		$O(\epsilon^{-(5-3\tau)/(2-6\tau)})$ Hessian-vector products.
	\end{itemize}
\end{theorem} 

\begin{remark}
	The initialization condition in Theorem~\ref{thm:bamm}
	can be satisfied by applying gradient descent to
	$\psi_{\sigma_0}(x^0,\cdot)$ with $x^0$ and $\sigma_0$ fixed.
	By strong convexity and smoothness, this requires
	$O(\sigma_0^{-1}\log(1/\sigma_0))$ LL gradient evaluations,
	with dependence on the initial distance and problem constants omitted.
	Since $k_0$ can be chosen independently of $K$,
	this one-time cost does not change the stated asymptotic
	oracle-complexity bounds.
\end{remark} 

We conclude this section by comparing the complexity bounds of the proposed methods with their corresponding baselines. Table~\ref{tab:com} summarizes the results. Outer-iteration complexities coincide with gradient/JVP complexities for all methods shown and are therefore omitted. Here, $\nabla_2f_{\sigma_k}^* :=\nabla_2f(x^k,y_{\sigma_k}^*(x^k))$. The GALET bounds are derived from the analysis in Appendix~F of \citet{xiao2023generalized}, as detailed in Appendix~\ref{app:galet-complexity}. The sl-BAMM bounds follow from Theorems~3.4--3.6 of \citet{liu2023averaged}.

\begin{table*}[t]
	\centering
	\renewcommand{\arraystretch}{1.25}
	\setlength{\tabcolsep}{4pt}
	\begin{small}
		\begin{tabular}{@{}ccccc@{}}
			\toprule
			Setting
			& Condition
			& Method
			& \makecell[c]{Gradient/JVP\\computations}
			& \makecell[c]{Hessian-vector\\products}\\
			\midrule
			
			\multirow{2}{*}{\makecell[c]{PL with positive\\spectral gap}}
			& \multirow{2}{*}{---}
			& GALET \citep{xiao2023generalized}
			& $O(\kappa^6\epsilon^{-1})$
			& $O(\kappa^8\epsilon^{-1}\log(\kappa/\epsilon))$\\
			\cmidrule(l){3-5}
			&
			& Newton-based GALET
			& $O(\kappa^4\epsilon^{-1})$
			& $O(\kappa^5\epsilon^{-1}\log(\kappa/\epsilon))$\\
			
			\midrule
			
			\multirow{6}{*}{\makecell[c]{$g(x,\cdot)$ convex,\\$f(x,\cdot)$ strongly convex}}
			& \multirow{2}{*}{\makecell[c]{$x^k,y^k,u^k$\\bounded}}
			& sl-BAMM \citep{liu2023averaged}
			& $O(\epsilon^{-11/(2-6\tau)})$
			& $O(\epsilon^{-11/(2-6\tau)})$\\
			\cmidrule(l){3-5}
			&
			& Newton-based BAMM
			& $O(\epsilon^{-4/(1-3\tau)})$
			& $O(\epsilon^{-(9-3\tau)/(2-6\tau)})$\\
			\cmidrule(l){2-5}
			
			& \multirow{2}{*}{\makecell[c]{$\nabla_2f_{\sigma_k}^*,u^k$\\bounded}}
			& sl-BAMM \citep{liu2023averaged}
			& $O(\epsilon^{-7/(2-6\tau)})$
			& $O(\epsilon^{-7/(2-6\tau)})$\\
			\cmidrule(l){3-5}
			&
			& Newton-based BAMM
			& $O(\epsilon^{-3/(1-3\tau)})$
			& $O(\epsilon^{-(7-3\tau)/(2-6\tau)})$\\
			\cmidrule(l){2-5}
			
			& \multirow{2}{*}{\makecell[c]{$u_{\sigma_k}^*,u^k$\\bounded}}
			& sl-BAMM \citep{liu2023averaged}
			& $O(\epsilon^{-5/(2-6\tau)})$
			& $O(\epsilon^{-5/(2-6\tau)})$\\
			\cmidrule(l){3-5}
			&
			& Newton-based BAMM
			& $O(\epsilon^{-2/(1-3\tau)})$
			& $O(\epsilon^{-(5-3\tau)/(2-6\tau)})$\\
			
			\bottomrule
		\end{tabular}
		\caption{Complexity comparison with the corresponding baselines beyond LL strong convexity. For BAMM, $O(\cdot)$ suppresses dependence on $\mu_f,L_{f,1},L_{f,2},L_{g,1},L_{g,2}$.}
		\label{tab:com}
	\end{small}
\end{table*}

\section{Numerical Experiments}
\label{experiments}

We evaluate the proposed methods on two toy examples, data hyper-cleaning, and GPT-2 data reweighting. 

\subsection{Toy Examples}

We first consider the positive spectral-gap setting:
\begin{align}
\begin{split}
f(x,y)
&=\sum_{i=1}^m\log\cosh(x_i-2)
+\sum_{i=1}^m\log\cosh([y_1]_i+3[y_2]_i-2[x]_i),\\
g(x,y)
&=\frac{\zeta}{2}\|y_1-y_2-x\|^2
+\sum_{i=1}^m\log\cosh([y_1]_i-[y_2]_i-[x]_i),
\end{split}
\label{eq:toy1}
\end{align}
where $x\in\mathbb{R}^m$, $y=(y_1,y_2)\in\mathbb{R}^{2m}$, and
$\cosh(t)=(e^t+e^{-t})/2$.
The LL objective is not strongly convex in $y$ and has solution set
$\{y:y_1=y_2+x\}$. Its positive Hessian eigenvalues are uniformly bounded below by $2\zeta$, so Assumption~\ref{asp:singularvalue} holds.

\begin{figure}[h]
	\centering
	\includegraphics[width=0.7\linewidth]{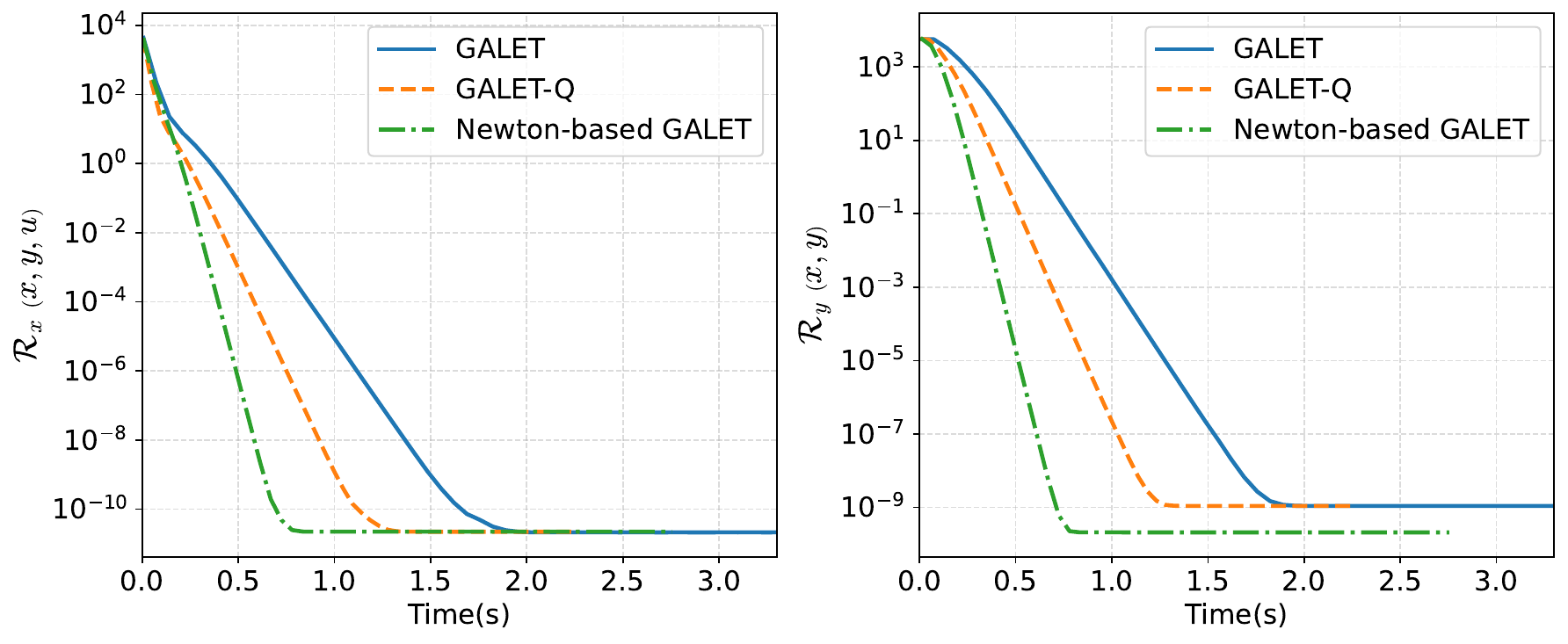}
	\caption{GALET, GALET-Q, and Newton-based GALET on \eqref{eq:toy1} with $\zeta=1$. \textbf{Left:} $\mathcal{R}_x$. \textbf{Right:} $\mathcal{R}_y$.}
	\label{fig:toy-galet1}
\end{figure}

\begin{figure}[H]
	\centering
	\includegraphics[width=0.7\linewidth]{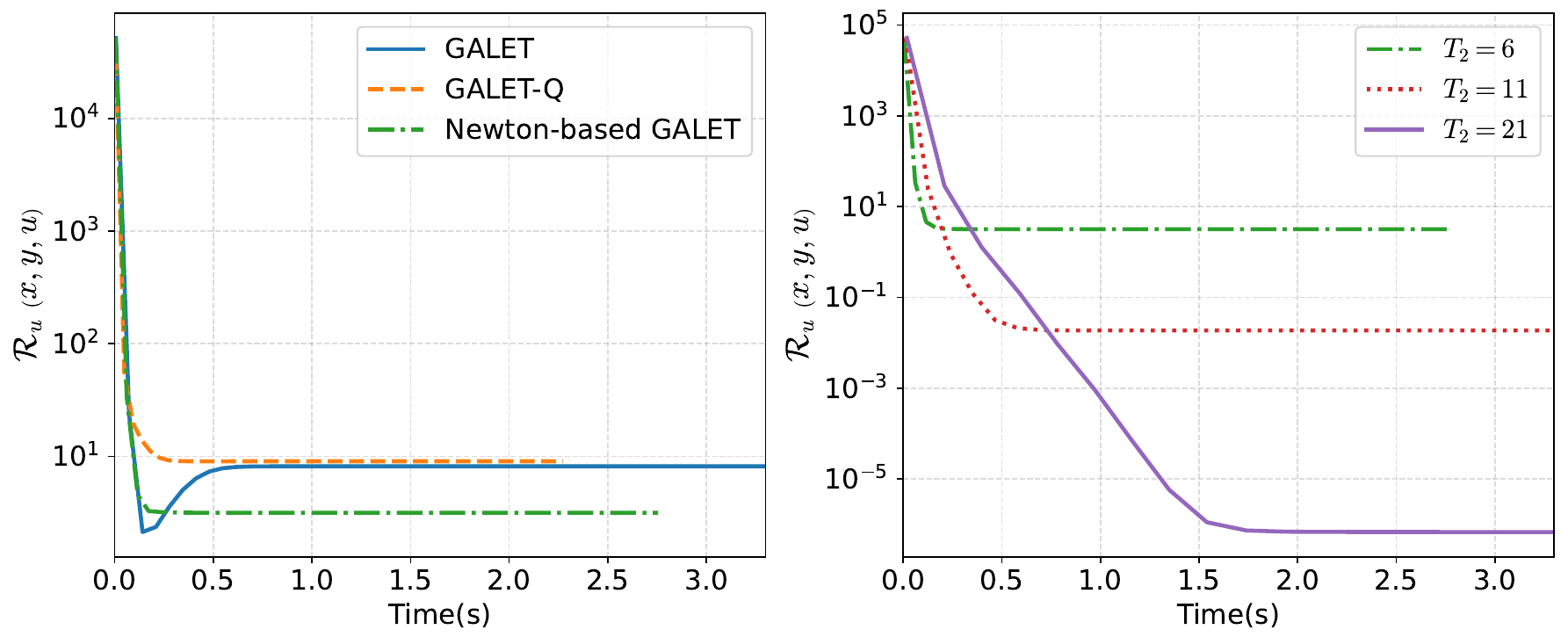}
	\caption{Results on \eqref{eq:toy1} with $\zeta=1$. \textbf{Left:} $\mathcal{R}_u$ for GALET, GALET-Q, and Newton-based GALET. \textbf{Right:} $\mathcal{R}_u$ for Newton-based GALET with $T_2=6,11,21$.}
	\label{fig:toy-galet2}
\end{figure}

Figures~\ref{fig:toy-galet1} and \ref{fig:toy-galet2} compare Newton-based GALET with GALET and GALET-Q, a variant that replaces Line~7 of Algorithm~1 in \cite{xiao2023generalized} by
\begin{equation*}
w^{k,t+1}
=
w^{k,t}
-\rho\left(
\nabla_2f(x^k,y^{k+1})
+\nabla_{22}^2g(x^k,y^{k+1})w^{k,t}
\right).
\end{equation*}
GALET-Q thus applies GD directly to the multiplier linear system rather than the squared-residual problem used in GALET, isolating the effect of this modification.

We set $N=1$ and $T=5$ for GALET and GALET-Q, and $T_1=2$, $T_2=6$ for Newton-based GALET. Since the GD subroutine starts from $p^0=0$, the first update $p^1=\gamma q$ is curvature-independent, so $T$ iterations contain $T-1$ curvature-dependent refinements. Although the analysis uses $T_1=\Theta(\kappa)$, we find that $T_1=2$ suffices in this example, consistent with the empirical observation in \cite{dong2025efficient}. We use additional iterations for the cold-start multiplier solve.

\begin{figure}[H]
	\centering
	\includegraphics[width=0.7\linewidth]{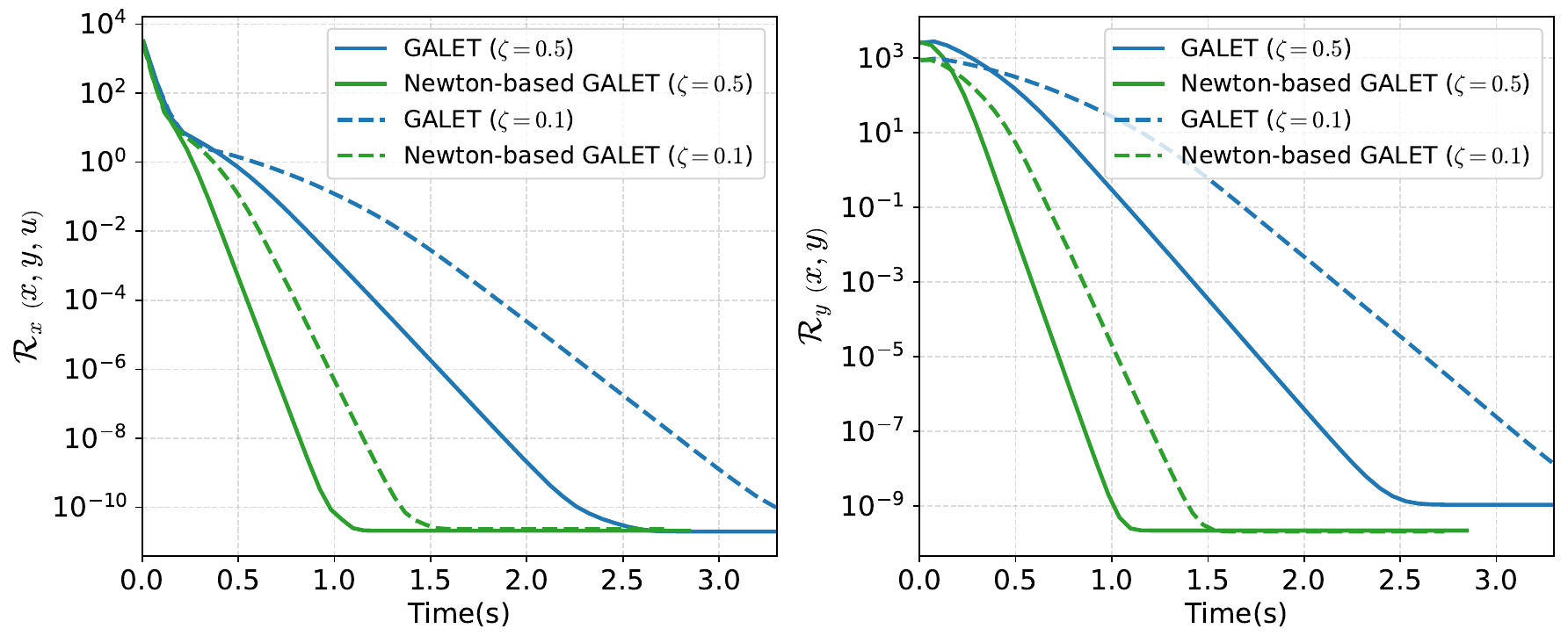}
	\caption{GALET and Newton-based GALET on \eqref{eq:toy1} with $\zeta=0.5,0.1$. \textbf{Left:} $\mathcal{R}_x$. \textbf{Right:} $\mathcal{R}_y$.}
	\label{fig:toy-galet3}
\end{figure}

Figure~\ref{fig:toy-galet1} shows faster residual reduction for GALET-Q than GALET and a further improvement for Newton-based GALET. Figure~\ref{fig:toy-galet2} shows that increasing $T_2$ reduces $\mathcal{R}_u$, consistent with Lemma~\ref{lem:u}. We also vary the condition number in \eqref{eq:toy1}, where $\kappa=\Theta(1/\zeta)$. Figure~\ref{fig:toy-galet3} shows that Newton-based GALET degrades less as $\kappa$ increases than GALET in this experiment.

For the second setting, with a convex LL objective and an UL objective strongly convex in the LL variable, consider
\begin{align}
\begin{split}
f(x,y)
&=\sum_{i=1}^m\log\cosh\bigl(x_i-[y_2]_i\bigr)
+\frac{1}{2}\|y_1-e\|^2+\frac{1}{2}\|y_2\|^2, \\
g(x,y)
&=\frac{1}{2}\|y_1\|^2-\langle x,y_1\rangle,
\end{split}
\label{eq:toy2}
\end{align}
where $x\in\mathbb{R}^m$, $y=(y_1,y_2)\in\mathbb{R}^{2m}$, and $e$ is the all-ones vector. The first term of $f$ is adapted from
\cite[Eq. (13)]{liu2023averaged} to ensure that
$\|\nabla_1f(x,y)\|$ is uniformly bounded. The LL solution set is
$\{y:y_1=x,\;y_2\in\mathbb{R}^m\}$, while $f(x,\cdot)$ is strongly convex.

\begin{figure}[H]
	\centering
	\includegraphics[width=0.7\linewidth]{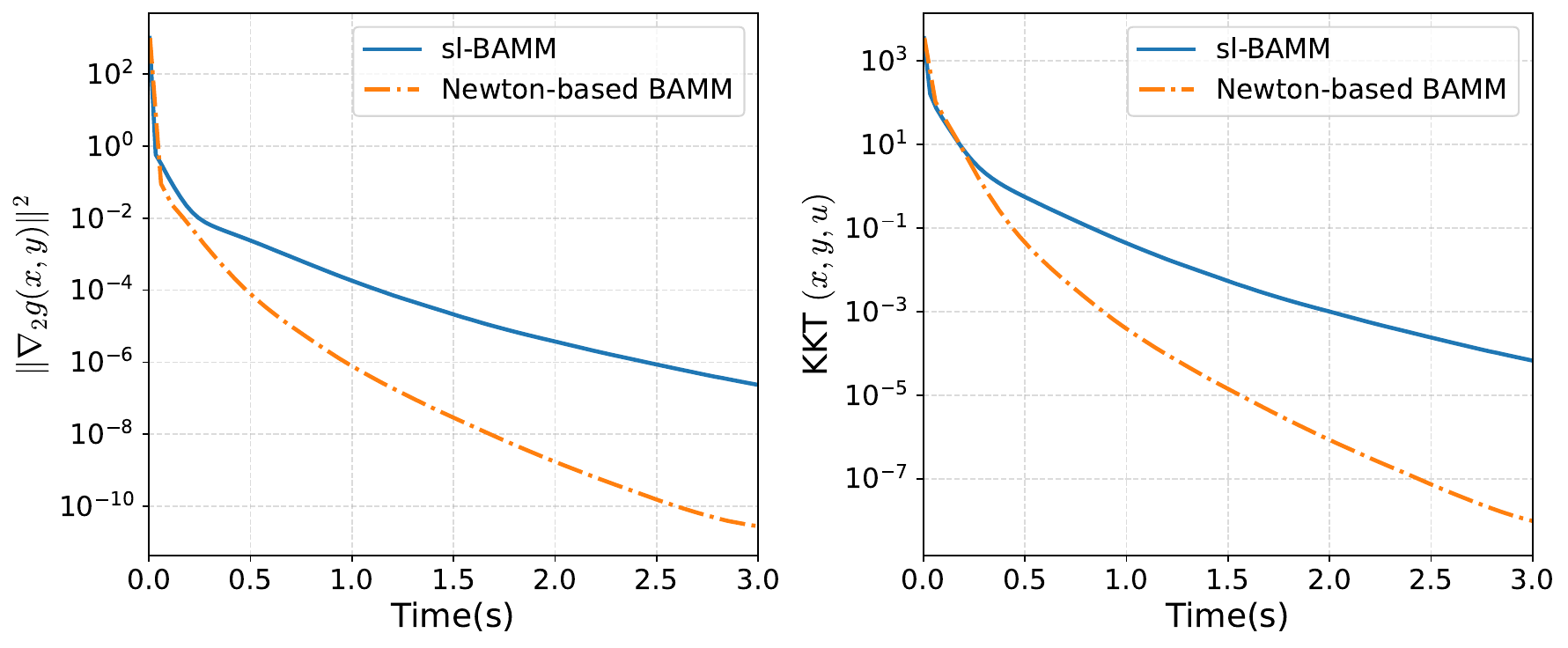}
	\caption{Newton-based BAMM and sl-BAMM on \eqref{eq:toy2}, with stepsizes chosen according to Case~1 of Theorem~\ref{thm:bamm}.}
	\label{fig:toy-bamm}
\end{figure}

Figure~\ref{fig:toy-bamm} compares Newton-based BAMM with sl-BAMM. We set $T_k=2$ for Newton-based BAMM. On this example, Newton-based BAMM shows faster reduction in both the LL and KKT residuals.

\subsection{Data Hyper-Cleaning}

We next consider data hyper-cleaning, a standard bilevel optimization benchmark. The UL variable $\lambda$ weights the training samples, while the LL variable $\omega$ represents the model parameters. We consider both linear and MLP models.

For the linear model, we use
\begin{equation*}
f(\lambda,\omega)=\frac{1}{|\mathcal D'|}\sum_{(x'_e,y'_e)\in\mathcal D'}\mathcal{L}(\omega x'_e,y'_e)+c_r\|\omega\|^2,\qquad
g(\lambda,\omega)=\frac{1}{|\mathcal D|}\sum_{(x_e,y_e)\in\mathcal D}\sigma(\lambda_e)\mathcal{L}(\omega x_e,y_e),
\end{equation*} 
where $\mathcal D$ and $\mathcal D'$ denote the training and validation sets, respectively, $\mathcal L$ is the cross-entropy loss, and $\sigma(\lambda_e)$ is the weight assigned to the $e$-th training sample. The LL objective is convex in $\omega$, while the quadratic regularizer with $c_r>0$ makes the UL objective strongly convex in $\omega$, matching the setting of Newton-based BAMM.

\begin{figure}[H]
	\centering
	\includegraphics[width=0.7\linewidth]{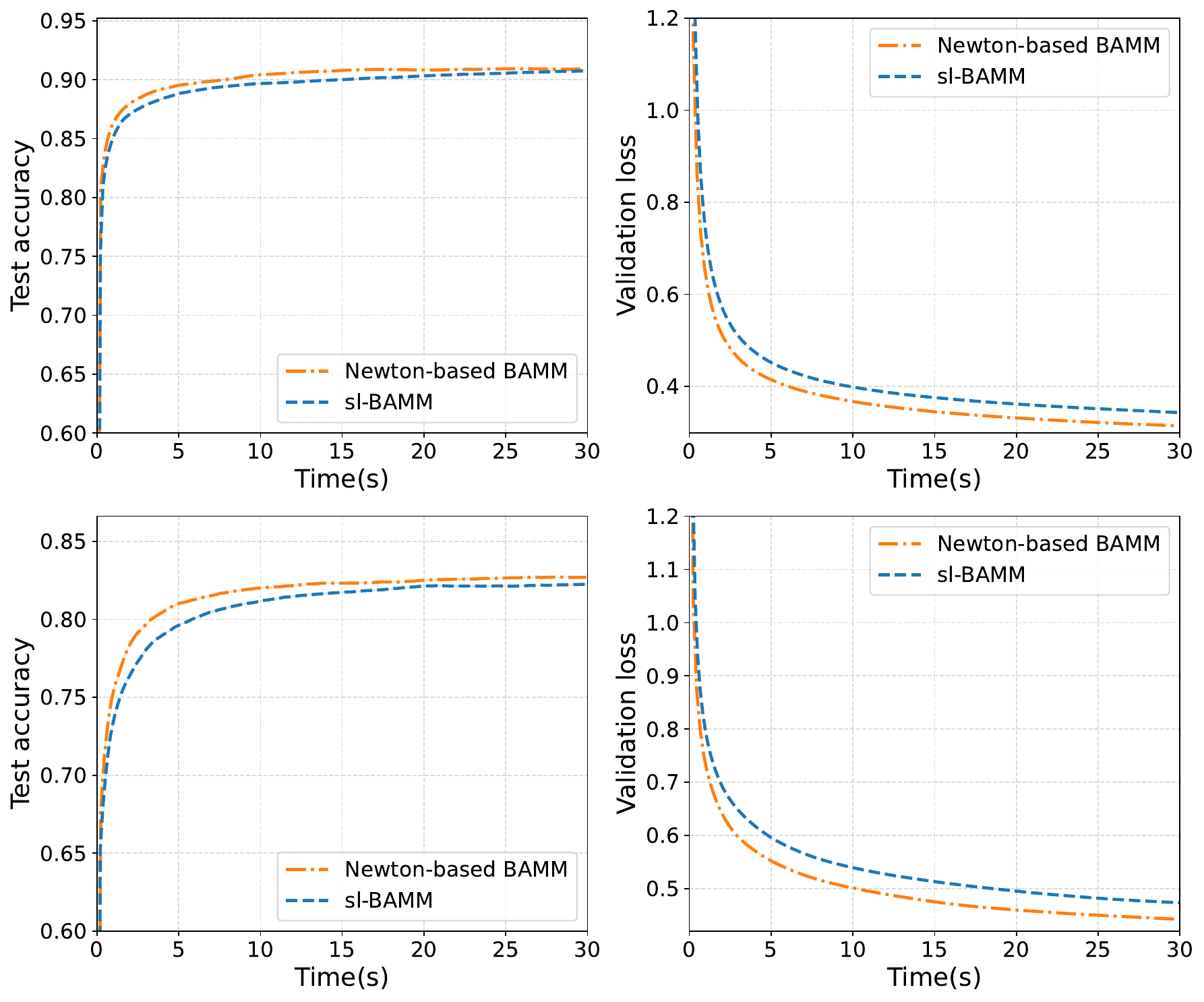}
	\caption{Data hyper-cleaning with linear models. Comparison of sl-BAMM and Newton-based BAMM on MNIST (upper) and FashionMNIST (lower).}
	\label{fig:hpo-bamm}
\end{figure}

We compare Newton-based BAMM with sl-BAMM on MNIST and FashionMNIST, using $T_k=2$ for Newton-based BAMM. Figure~\ref{fig:hpo-bamm} shows faster reduction in the reported metrics for Newton-based BAMM on both datasets.

\begin{figure}[h]
	\centering
	\includegraphics[width=0.7\linewidth]{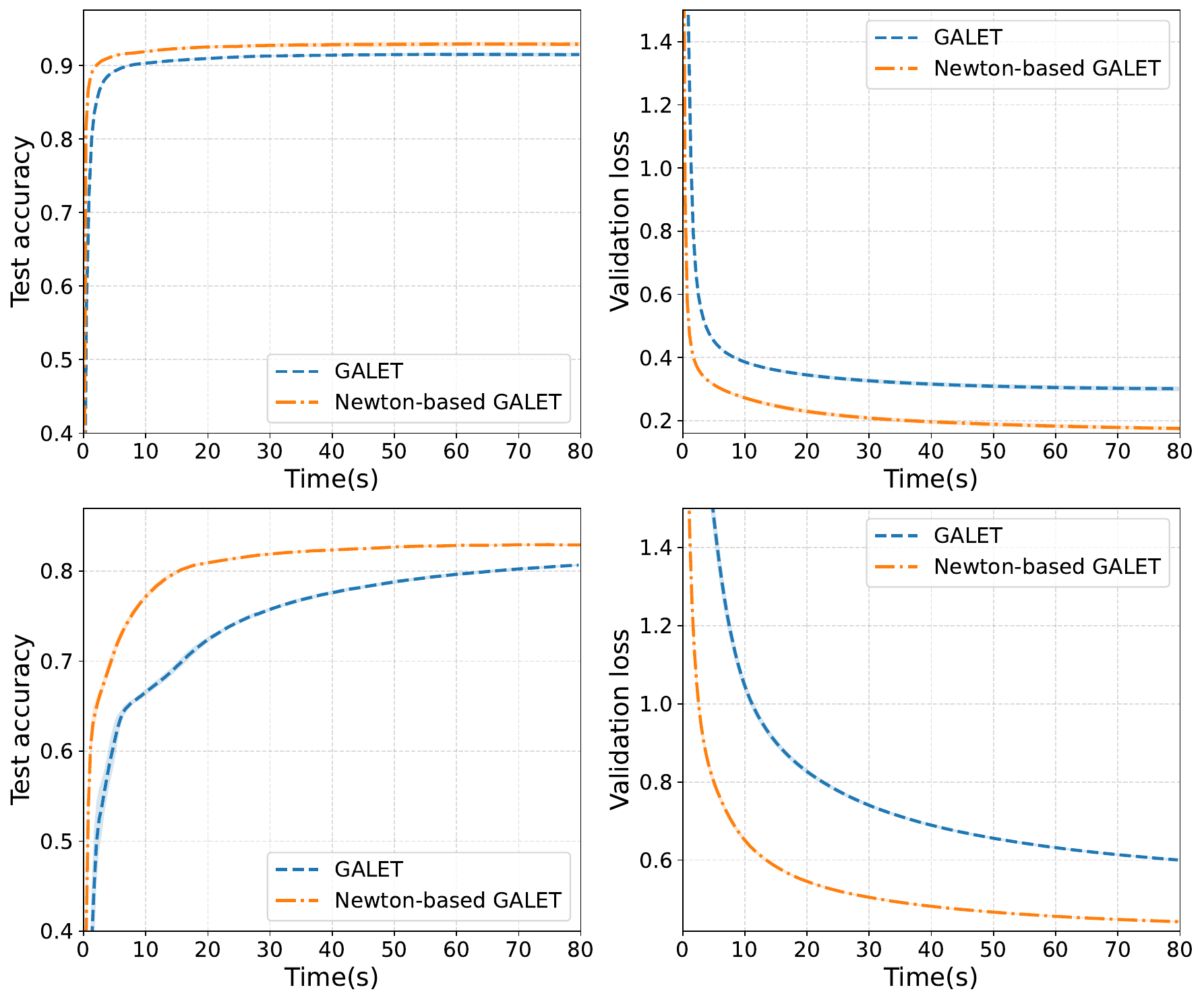}
	\caption{Data hyper-cleaning with two-layer MLPs. Comparison of GALET and Newton-based GALET on MNIST (upper) and FashionMNIST (lower).}
	\label{fig:hpo-galet}
\end{figure}

For the MLP model, we consider
\begin{equation*}
f(\lambda,\omega)=\frac{1}{|\mathcal D'|}\sum_{(x'_e,y'_e)\in\mathcal D'}\mathcal{L}(N(x'_e,\omega),y'_e),\qquad
g(\lambda,\omega)=\frac{1}{|\mathcal D|}\sum_{(x_e,y_e)\in\mathcal D}\sigma(\lambda_e)\mathcal{L}(N(x_e,\omega),y_e),
\end{equation*}
where $N(\cdot,\omega)$ is a two-layer MLP. Unlike the linear setting, neither the convexity of the LL objective nor the strong convexity of the UL objective in $\omega$ is guaranteed. We therefore compare Newton-based GALET with GALET. Figure~\ref{fig:hpo-galet} shows that Newton-based GALET achieves higher accuracy within the reported running time on both datasets.

\vspace{-0.8cm}


\subsection{Data Denoising with GPT-2}

We further evaluate Newton-based GALET on the data-denoising task of
\cite{pan2025scalebio}, using GPT-2 with 124M parameters. The validation set
$\mathcal{D}_{\mathrm{val}}$ contains 1,000 clean Alpaca samples
\citep{taori2023stanford}, while the training set consists of 1,000 clean samples
$\mathcal{D}_1$ and 9,000 corrupted samples $\mathcal{D}_2$, whose responses are replaced by ``.''. The goal is to assign a large weight to $\mathcal{D}_1$ and a small weight to $\mathcal{D}_2$.

We formulate the task as
\begin{equation*}
f(\lambda,w)=\frac{1}{|\mathcal{D}_{\mathrm{val}}|}\sum_{z\in\mathcal{D}_{\mathrm{val}}}\ell(w;z),\qquad
g(\lambda,w)=\sum_{i=1}^{2}\frac{p_i(\lambda)}{|\mathcal{D}_i|}\sum_{z\in\mathcal{D}_i}\ell(w;z),
\end{equation*} 
where $w$ denotes the GPT-2 parameters, $\ell(w;z)$ is the causal language modeling loss, and
$p_i(\lambda):=\exp(\lambda_i)/\sum_{j=1}^{2}\exp(\lambda_j)$
is the sampling weight of $\mathcal{D}_i$.

\begin{figure}[!htbp]
	\centering
	\includegraphics[width=0.5\linewidth]{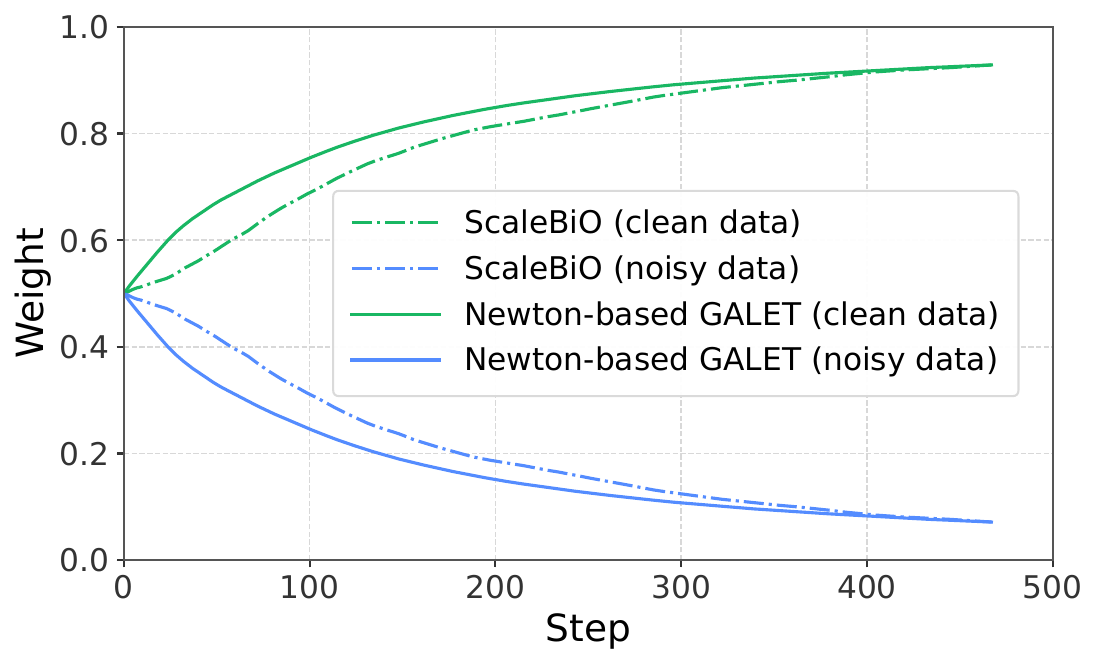}
	\caption{Clean and corrupted data weights for ScaleBiO and Newton-based GALET, averaged over five random seeds.}
	\label{fig:denoising_weight}
\end{figure}

We compare Newton-based GALET with ScaleBiO \citep{pan2025scalebio}, a single-loop first-order bilevel method. For Newton-based GALET, we set $T_1=T_2=2$; both methods use AdamW for all variable updates. Figure~\ref{fig:denoising_weight} shows that both methods increase the clean-source weight and decrease the corrupted-source weight, with Newton-based GALET adjusting the weights faster initially. After three epochs, Table~\ref{tab:denoising_results} shows similar source weights but lower validation and test losses for Newton-based GALET, together with higher peak GPU memory due to the additional second-order computations.

\begin{table}[H]
	\centering
	\small
	\setlength{\tabcolsep}{3.5pt}
	\begin{tabular}{cccccc}
		\toprule
		Method
		& Clean weight
		& Corrupted weight
		& Validation loss
		& Test loss
		& Peak memory\\
		\midrule
		ScaleBiO
		& 0.928
		& 0.072
		& 2.603
		& 2.555
		& 9.984 GiB\\
		\midrule
		\makecell[c]{Newton-based\\GALET}
		& 0.929
		& 0.071
		& 2.439
		& 2.398
		& 18.509 GiB\\
		\bottomrule
	\end{tabular}
	\caption{GPT-2 data denoising after three training epochs, averaged over five random seeds.}
	\label{tab:denoising_results}
\end{table} 


\section{A Generalized Newton-Based Framework}
\label{sec:discussion}

The proposed methods share a common structure: the LL and multiplier updates exploit curvature from the same LL model. Newton-based GALET uses $g$ with potentially singular curvature, whereas Newton-based BAMM uses the strongly convex aggregation model $\psi_{\sigma_k}$. Related structures also appear in Newton and quasi-Newton bilevel methods
\citep{ramzi2022shine,fang2025qnbo,dong2025efficient}. We capture these common features in a generalized Newton-based framework for bilevel optimization (GNBO).

\paragraph{Lower-level model}
At iteration $k$, GNBO selects an LL model $h_k(x,y)$ with residuals
\begin{align*}
r_y(x,y)&=\nabla_2h_k(x,y),\\
r_u(x,y,u)&=\nabla_2f(x,y)-\nabla_{22}^2h_k(x,y)u,\\
r_x(x,y,u)&=\nabla_1f(x,y)-\nabla_{12}^2h_k(x,y)u.
\end{align*}
Newton-based GALET takes $h_k=g$, while Newton-based BAMM takes
$h_k=\psi_{\sigma_k}$; in both cases, $h_k$ supplies the curvature shared by the LL and multiplier updates.

\paragraph{Lower-level update}
A Newton-type LL update targets
$\nabla_{22}^2h_k(x^k,y)v=\nabla_2h_k(x^k,y)$.
To allow one or multiple LL corrections per outer iteration, let
$y^{0,k}=y^k,\ldots,y^{S,k}=y^{k+1}$ and consider
\begin{equation}
B_v^{s,k}v=q_v^{s,k},
\label{eq:linear-system-v}
\end{equation}
where $B_v^{s,k}$ approximates the curvature and $q_v^{s,k}$ is derived from the LL residual. A subroutine $\mathcal{A}_v$ approximately solves
\eqref{eq:linear-system-v}, followed by
$y^{s+1,k}=y^{s,k}-\beta_{s,k}v^{s,k}$.
Here, $S$ counts LL corrections and $T_1^k$ denotes the internal solver effort.

The multiplier and UL updates may use either the pre- or post-update LL point:
$\tilde y^k=\delta_y y^k+(1-\delta_y)y^{k+1}$,
where $\delta_y\in\{0,1\}$. Thus, $\delta_y=0$ uses $y^{k+1}$, while
$\delta_y=1$ uses $y^k$, as in Newton-based GALET, Newton-based BAMM, and NBO-GD.

\paragraph{Multiplier update}
At $\tilde y^k$, the ideal multiplier system is
$\nabla_{22}^2h_k(x^k,\tilde y^k)u=\nabla_2f(x^k,\tilde y^k)$.
GNBO represents its treatment by
\begin{equation}
B_u^k(\delta_u u^k-w)=q_u^k,
\label{eq:linear-system-w}
\end{equation}
where $B_u^k$ and $q_u^k$ specify the curvature approximation and right-hand side. The indicator $\delta_u\in\{0,1\}$ distinguishes cold and warm starts:
$\delta_u=0$ gives $u^{k+1}=-w^k$, as in Newton-based GALET, while
$\delta_u=1$ gives $u^{k+1}=u^k-w^k$, as in Newton-based BAMM and NBO-GD.
A subroutine $\mathcal{A}_w$ approximately solves \eqref{eq:linear-system-w}, with $T_2^k$ denoting the internal solver effort.

\paragraph{Upper-level update}
Set $\tilde u^k=\delta_u u^k+(1-\delta_u)u^{k+1}$ and
\begin{align*}
d_x^k
=
\nabla_1f(x^k,\tilde y^k)
-\nabla_{12}^2h_k(x^k,\tilde y^k)\tilde u^k.
\end{align*}
Then $x^{k+1}=x^k-\alpha_kd_x^k$.
Algorithm~\ref{alg:GNBO} summarizes the framework. 

\begin{algorithm}[h]
	\caption{Generalized Newton-Based Bilevel Optimization (GNBO)}
	\label{alg:GNBO}
	\begin{algorithmic}[1]
		\STATE \textbf{Input:} $x^0,y^0,u^0$; $K,S,\{T_1^k,T_2^k\}$; step sizes $\{\alpha_k,\beta_{s,k}\}$; $\delta_y,\delta_u\in\{0,1\}$.
		\FOR{$k=0,\ldots,K-1$}
		\STATE Choose an LL model $h_k$ and set $y^{0,k}=y^k$. \COMMENT{Lower-Level Update}
		\FOR{$s=0,\ldots,S-1$}
		\STATE $v^{s,k}=\mathcal{A}_v(h_k;x^k,y^{s,k};T_1^k)$ based on \eqref{eq:linear-system-v}.
		\STATE $y^{s+1,k}=y^{s,k}-\beta_{s,k}v^{s,k}$.
		\ENDFOR
		\STATE Set $y^{k+1}=y^{S,k}$ and $\tilde y^k=\delta_y y^k+(1-\delta_y)y^{k+1}$.
		\STATE $w^k=\mathcal{A}_w(h_k;x^k,\tilde y^k,u^k;T_2^k)$ based on \eqref{eq:linear-system-w}. \COMMENT{Multiplier Update}
		\STATE Set $u^{k+1}=\delta_u u^k-w^k$ and $\tilde u^k=\delta_u u^k+(1-\delta_u)u^{k+1}$.
		\STATE Set $d_x^k=\nabla_1f(x^k,\tilde y^k)-\nabla_{12}^2h_k(x^k,\tilde y^k)\tilde u^k$. \COMMENT{Upper-Level Update}
		\STATE Set $x^{k+1}=x^k-\alpha_kd_x^k$.
		\ENDFOR
	\end{algorithmic}
\end{algorithm} 

The number of LL corrections $S$ is distinct from the internal solver effort: $S=1$ means one LL update per outer iteration, even if computing its direction requires multiple solver iterations. Table~\ref{tab:GNBO-instances} summarizes representative instances of GNBO. 

\begin{table}[t!]
	\centering
	\setlength{\tabcolsep}{2.5pt}
	\renewcommand{\arraystretch}{1.08}
	\begin{tabular}{@{}ccccc@{}}
		\toprule
		Method & LL model & LL corrections & Solver / curvature & $(\delta_y,\delta_u)$\\
		\midrule
		\makecell{SHINE\\ \citep{ramzi2022shine}}
		& $g$ & Multiple & Quasi-Newton, reused & $(0,0)$\\
		\makecell{qNBO\\ \citep{fang2025qnbo}}
		& $g$ & Multiple & Separate quasi-Newton & $(0,0)$\\
		\makecell{NBO-GD \\\citep{dong2025efficient}}
		& $g$ & One & \makecell{$T$ GD steps,\\ common Hessian} & $(1,1)$\\
		\midrule
		Newton-based GALET
		& $g$ & One & \makecell{$T_1/T_2$ GD steps,\\ common Hessian} & $(1,0)$\\
		Newton-based BAMM
		& $\psi_{\sigma_k}$ & One & \makecell{$T_k$ GD steps,\\ common Hessian} & $(1,1)$\\
		\bottomrule
	\end{tabular}
	\caption{Representative instances of GNBO.}
	\label{tab:GNBO-instances}
\end{table}


\section{Concluding Remarks}
\label{conclude}

We developed Newton-based methods for bilevel optimization beyond LL strong convexity and captured their shared-curvature structure through GNBO. For the settings studied, exploiting the same curvature in the LL and multiplier updates leads to improved complexity bounds over the corresponding baselines.

Several directions remain open. Hessian-vector products introduce additional memory and computational costs, particularly for large-scale neural networks. Improving scalability through parameter-efficient strategies such as LoRA \citep{hu2022lora} and LISA \citep{pan2024lisa}, or limited-memory quasi-Newton approximations \citep{liu1989limited}, is therefore of interest. Moreover, our methods use Newton-type LL and multiplier updates but retain a first-order UL update. Extending this approach to Newton-type UL updates under suitable higher-order smoothness conditions, as explored in recent work \citep{bi2026second,yang2026second}, is another direction for future study.

\newpage

\appendix

\section{Experimental Details}

All experiments are conducted on a machine with an Intel Xeon Gold 5218R CPU (2.10 GHz) and an NVIDIA H100 GPU (80 GB). The toy and data hyper-cleaning experiments are repeated over 10 runs, while the GPT-2 experiments use five random seeds.

\paragraph{Toy examples}
For GALET, GALET-Q, and Newton-based GALET, we set
$\alpha_k=\eta_k=0.1$ and $\gamma_k=0.01$.
For sl-BAMM, we use $p=37/440$ and $\tau=1/40$;
for Newton-based BAMM, we use $p=37/320$ and $\tau=1/40$.

\paragraph{Data hyper-cleaning}
For GALET and Newton-based GALET, we jointly tune
$(\alpha_k,\beta_k,\eta_k)$ by grid search, with
$\alpha_k\in\{10^3,10^4,10^5\}$ and
$\gamma_k,\eta_k\in\{10^{-3/2},10^{-1},10^{-1/2}\}$
for MNIST or
$\{10^{-5/2},10^{-2},10^{-3/2}\}$
for FashionMNIST.
For sl-BAMM and Newton-based BAMM, we use
$(p,\tau)=(37/440,1/40)$ and $(37/320,1/40)$, respectively.

\paragraph{Data denoising with GPT-2}
ScaleBiO uses AdamW with a constant learning rate of $10^{-5}$
for both model-variable updates. Newton-based GALET uses AdamW
with the same learning rate for both linear-system solves.
In both methods, the source-weight logits are updated using
AdamW with a constant learning rate of $10^{-2}$.
GPT-2 parameters are stored in BF16, with LISA disabled. 


\section{Proofs of Main Results and Auxiliary Lemmas}
\label{apd:proof}

This appendix provides detailed proofs of the main convergence results, together with the statements and proofs of the auxiliary lemmas used in the analysis.

\subsection{Convergence Analysis of Newton-Based GALET} 

\begin{proof}[Proof of Lemma~\ref{lem:range}]
	Fix $x$. By Weyl's inequality and the $L_{g,2}$-Lipschitz continuity of
	$\nabla_{22}^2g(x,\cdot)$, its eigenvalues vary continuously with $y$.
	Assumption~\ref{asp:singularvalue} prevents any nonzero eigenvalue from
	crossing zero, so the inertia and rank of $\nabla_{22}^2g(x,y)$ are
	independent of $y$. At any $y^*\in Y^*(x)$, the second-order necessary
	condition gives $\nabla_{22}^2g(x,y^*)\succeq0$. Hence
	$\nabla_{22}^2g(x,y)\succeq0$ for all $y$, and $g(x,\cdot)$ is convex.
	Denote the common rank by $r(x)$.
	
	Apply \cite[Lemma~2]{hartman1959spherical} to
	$p(y):=\nabla_2g(x,y)$, whose Jacobian has constant rank $r(x)$.
	Since $p(y^*)=0$, its zero set is locally an affine plane
	$y^*+L(x)$ of dimension $n-r(x)$. Convexity of
	$Y^*(x)=\{y:p(y)=0\}$ extends this local characterization globally:
	$Y^*(x)=y^*+L(x)$. Moreover, differentiating
	$p(y^*+td)=0$ at $t=0$ for $d\in L(x)$ gives
	$L(x)\subseteq\ker(\nabla_{22}^2g(x,y^*))$. Since both spaces have
	dimension $n-r(x)$,
	\begin{equation}
	Y^*(x)=y^*+K(x),
	\qquad
	K(x):=\ker(\nabla_{22}^2g(x,y^*)).
	\label{eq:Y-affine}
	\end{equation}
	
	Let $S(x):=K(x)^\perp$. For $d\in K(x)$, convexity and
	\eqref{eq:Y-affine} imply
	$\langle\nabla_2g(x,y),y-y^*-td\rangle\geq0$ for every $t\in\mathbb R$.
	Hence $\nabla_2g(x,y)\perp K(x)$, so
	$\nabla_2g(x,y)\in S(x)$. Differentiating this relation with respect to
	$y$ gives
	$\range{\nabla_{22}^2g(x,y)}\subseteq S(x)$.
	Both spaces have dimension $r(x)$; therefore
	\begin{equation}
	\range{\nabla_{22}^2g(x,y)}=S(x),
	\qquad
	\nabla_2g(x,y)\in S(x).
	\label{eq:H-range}
	\end{equation}
	
	Let $\hat y$ be the projection of $y$ onto $Y^*(x)$. By
	\eqref{eq:Y-affine}, $y-\hat y\in S(x)$. Since the positive eigenvalues
	of the Hessian are at least $\mu$, \eqref{eq:H-range} gives
	\begin{equation*}
	\langle\nabla_2g(x,y),y-\hat y\rangle
	=\int_0^1 (y-\hat y)^\top\nabla_{22}^2g\bigl(x,\hat y+s(y-\hat y)\bigr)(y-\hat y)\,\mathrm{d}s
	\geq\mu\|y-\hat y\|^2.
	\end{equation*} 
	Thus, by Cauchy--Schwarz,
	$\|\nabla_2g(x,y)\|
	\geq\mu\|y-\hat y\|
	=\mu\dist{y,Y^*(x)}$.
	
	Finally, suppose $\nabla_2g(\cdot,y)$ is $L_{g,1}$-Lipschitz.
	For $d\in K(x)$, \eqref{eq:Y-affine} and the error bound give, for any
	$x'$ and $t\in\mathbb R$,
	\begin{equation*}
	\dist{y^*(x)+td,Y^*(x')}
	\leq\frac{L_{g,1}}{\mu}\|x'-x\|.
	\end{equation*}
	Since $Y^*(x')=y^*(x')+K(x')$, dividing by $t>0$ and letting
	$t\to\infty$ yields $d\in K(x')$. Interchanging $x$ and $x'$ shows that
	$K(x)=K(x')=:K$, so $S:=K^\perp$ is independent of $x$.
	Since $\nabla_2g(x,y)\in S$ for all $x,y$, differentiating with respect
	to $x$ gives
	$\range{\nabla_{21}^2g(x,y)}\subseteq S$.
	By \eqref{eq:H-range},
	$S=\range{\nabla_{22}^2g(x,y)}$, and hence $\range{\nabla_{21}^2g(x,y)}
	\subseteq
	\range{\nabla_{22}^2g(x,y)}.$
\end{proof} 

Following in part \cite{xiao2023generalized}, define
\[
F(x,y;u):=f(x,y)-u^\top\nabla_2g(x,y),
\]
with $u$ treated as fixed, and
$u^*(x,y):=[\Hxy]^\dagger\nabla_2f(x,y)$.
For the subsequent analysis, let
\begin{equation}
L_F:=L_{f,1}+\frac{L_{f,0}L_{g,2}}{\mu},
\qquad
L_{u^*}:=
\frac{(1+\sqrt{5})L_{f,0}L_{g,2}}{2\mu^2}
+\frac{L_{f,1}}{\mu}.
\label{eq:constant}
\end{equation} 

\begin{lemma}\label{lem:error-u}
	Suppose Assumptions~\ref{asp:singularvalue} and \ref{asp:smooth1} hold. Then
	\begin{equation*}
	\|u^*(x,y)-u^*(x',y')\|
	\leq
	L_{u^*}\bigl(\|x-x'\|+\|y-y'\|\bigr).
	\end{equation*}
\end{lemma}

\begin{proof}
	By Lemma~\ref{lem:range}, $\nabla_{22}^2g(x,y)$ has a common range, and hence constant rank, over $(x,y)$. The result follows from the pseudoinverse perturbation bound
	\cite[Theorem~3.3]{stewart1977perturbation} and the Lipschitz continuity of
	$\nabla_{22}^2g$ and $\nabla_2f$; see also
	\cite[Lemma~9]{xiao2023generalized}.
\end{proof} 

\begin{lemma}\label{lem:F}
	Suppose Assumptions~\ref{asp:singularvalue} and \ref{asp:smooth1} hold. Then Newton-based GALET satisfies
	\begin{align*}
	&F(x^{k+1},y^{k+1};u^*(x^{k+1},y^{k+1}))
	-F(x^k,y^k;u^*(x^k,y^k))\\
	\leq{}&
	\left(
	-\alpha_k+\frac{\alpha_k^2}{2}
	+\frac{L_F\alpha_k^2}{2}
	+\mu^2L_{u^*}^2\alpha_k^2
	+\frac{3L_{g,1}^2\alpha_k^2}{2\mu^2}
	\right)\|d_x^k\|^2
	+\frac{3}{2\mu^2}\|\nabla_2g(x^k,y^k)\|^2\\
	&+
	\left(
	\frac{L_F}{2}
	+\mu^2L_{u^*}^2
	+\frac{3L_{g,1}^2}{2\mu^2}
	\right)\|v^k\|^2
	+\frac{L_{g,1}^2}{2\mu^2}
	\|\Hk u^{k+1}-\Hk u^*(x^k,y^k)\|^2.
	\end{align*}
\end{lemma}

\begin{proof}
	Let $g_k:=\nabla_2g(x^k,y^k)$ and $u_k^*:=u^*(x^k,y^k)$. For fixed $u$, $F(\cdot,\cdot;u)$ has an
	$(L_{f,1}+L_{g,2}\|u\|)$-Lipschitz continuous gradient. Since
	$\|u_k^*\|\leq L_{f,0}/\mu$, $F(\cdot,\cdot;u_k^*)$ has an $L_F$-Lipschitz continuous gradient.
	Decompose the potential difference as $\Delta_1+\Delta_2$, where
	\begin{align*}
	\Delta_1
	&:=F(x^{k+1},y^{k+1};u_k^*)-F(x^k,y^k;u_k^*),\\
	\Delta_2
	&:=F(x^{k+1},y^{k+1};u^*(x^{k+1},y^{k+1}))
	-F(x^{k+1},y^{k+1};u_k^*).
	\end{align*}
	By $L_F$-smoothness,
	\begin{align}
	\Delta_1
	\leq{}&
	\langle\nabla_1F(x^k,y^k;u_k^*),x^{k+1}-x^k\rangle
	+\langle\nabla_2F(x^k,y^k;u_k^*),y^{k+1}-y^k\rangle \nonumber\\
	&+\frac{L_F}{2}
	\left(\|x^{k+1}-x^k\|^2+\|y^{k+1}-y^k\|^2\right).
	\label{eq:F-delta1}
	\end{align}
	
	By Lemma \ref{lem:range},
	$\range{\nabla_{21}^2g(x^k,y^k)}\subseteq\range\Hk$.
	Since $\Hk$ is symmetric, this implies
	$\nabla_{12}^2g(x^k,y^k)
	=\nabla_{12}^2g(x^k,y^k)[\Hk]^\dagger\Hk$.
	Using this identity and $x^{k+1}-x^k=-\alpha_kd_x^k$,
	\begin{align*}
	&\langle\nabla_1F(x^k,y^k;u_k^*),x^{k+1}-x^k\rangle
	= -\alpha_k\|d_x^k\|^2\\
	&\qquad+\left\langle
	\nabla_{12}^2g(x^k,y^k)[\Hk]^\dagger
	\bigl(\Hk u^{k+1}-\Hk u_k^*\bigr),
	x^{k+1}-x^k
	\right\rangle\\
	&\leq{}
	-\left(\alpha_k-\frac{\alpha_k^2}{2}\right)\|d_x^k\|^2
	+\frac{L_{g,1}^2}{2\mu^2}
	\|\Hk u^{k+1}-\Hk u_k^*\|^2,
	\end{align*}
	where we used
	$\|\nabla_{12}^2g(x^k,y^k)\|_{\operatorname{op}}\leq L_{g,1}$ and
	$\|[\Hk]^\dagger\|_{\operatorname{op}}\leq1/\mu$.
	For the $y$-term, Algorithm~\ref{alg:PLGD} gives
	\begin{equation}
	v^k
	=
	\gamma_k\sum_{i=0}^{T_1-1}
	(I-\gamma_k\Hk)^i g_k
	\in\range\Hk.
	\label{eq:v}
	\end{equation}
	Let $P_k:=\Hk[\Hk]^\dagger$ be the orthogonal projector onto $\range\Hk$. Since
	$\Hk u_k^*=P_k\nabla_2f(x^k,y^k)$,
	\begin{equation*}
	\nabla_2F(x^k,y^k;u_k^*)
	=
	(I-P_k)\nabla_2f(x^k,y^k)
	\in(\range\Hk)^\perp.
	\end{equation*}
	As $y^{k+1}-y^k=-v^k\in\range\Hk$, the corresponding inner product in \eqref{eq:F-delta1} vanishes. Hence
	\begin{align}
	\Delta_1
	\leq &
	\left(-\alpha_k+\frac{\alpha_k^2}{2}
	+\frac{L_F\alpha_k^2}{2}\right)\|d_x^k\|^2
	+\frac{L_F}{2}\|v^k\|^2 \nonumber\\
	&+\frac{L_{g,1}^2}{2\mu^2}
	\|\Hk u^{k+1}-\Hk u_k^*\|^2.
	\label{eq:delta1-bound}
	\end{align}
	
	For $\Delta_2$, linearity of $F$ in $u$ and Young's inequality give
	\begin{align*}
	\Delta_2
	&=
	\left\langle
	u_k^*-u^*(x^{k+1},y^{k+1}),
	\nabla_2g(x^{k+1},y^{k+1})
	\right\rangle\\
	&\leq
	\frac{\mu^2}{2}
	\|u^*(x^{k+1},y^{k+1})-u_k^*\|^2
	+\frac{1}{2\mu^2}
	\|\nabla_2g(x^{k+1},y^{k+1})\|^2.
	\end{align*}
	By Lemma~\ref{lem:error-u} and the updates of $x$ and $y$,
	\begin{equation*}
	\|u^*(x^{k+1},y^{k+1})-u_k^*\|^2
	\leq
	2L_{u^*}^2
	\left(\alpha_k^2\|d_x^k\|^2+\|v^k\|^2\right).
	\end{equation*}
	Moreover, the $L_{g,1}$-Lipschitz continuity of $\nabla_2g$ gives
	\begin{equation*}
	\|\nabla_2g(x^{k+1},y^{k+1})\|^2
	\leq
	3\|g_k\|^2
	+3L_{g,1}^2\alpha_k^2\|d_x^k\|^2
	+3L_{g,1}^2\|v^k\|^2.
	\end{equation*}
	Therefore,
	\begin{align}
	\Delta_2
	\leq
	\left(
	\mu^2L_{u^*}^2\alpha_k^2
	+\frac{3L_{g,1}^2\alpha_k^2}{2\mu^2}
	\right)\|d_x^k\|^2
	+\frac{3}{2\mu^2}\|g_k\|^2 
	+
	\left(
	\mu^2L_{u^*}^2
	+\frac{3L_{g,1}^2}{2\mu^2}
	\right)\|v^k\|^2.
	\label{eq:delta2-bound}
	\end{align}
	Combining \eqref{eq:delta1-bound} and \eqref{eq:delta2-bound} proves the result.
\end{proof}

\begin{lemma}\label{lem:y}
	Suppose Assumptions~\ref{asp:singularvalue} and \ref{asp:smooth1} hold. In Newton-based GALET, let
	$0 \leq \gamma_k\leq1/L_{g,1}$ and
	$T_1\geq\log(1/4)/\log(1-\mu\gamma_k)$.
	If
	$\|\nabla_2g(x^k,y^k)\|\leq\mu^2/(2L_{g,2})$,
	then
	\begin{equation*}
	\|\nabla_2g(x^{k+1},y^{k+1})\|
	\leq
	\frac{1}{2}\|\nabla_2g(x^k,y^k)\|
	+L_{g,1}\|x^{k+1}-x^k\|.
	\end{equation*}
\end{lemma}

\begin{proof}
	Let $g_k:=\nabla_2g(x^k,y^k)$. Since
	$y^{k+1}=y^k-v^k$ and $\nabla_{22}^2g(x^k,\cdot)$ is $L_{g,2}$-Lipschitz continuous, Taylor's theorem gives
	\begin{equation}
	\|\nabla_2g(x^k,y^{k+1})-g_k+\Hk v^k\|
	\leq
	\frac{L_{g,2}}{2}\|v^k\|^2.\nonumber
	\end{equation}
	Define $r_k:=g_k-\Hk v^k$. Then
	\begin{equation}
	\|\nabla_2g(x^k,y^{k+1})\|
	\leq
	\frac{L_{g,2}}{2}\|v^k\|^2+\|r_k\|.
	\label{eq:y-descent0}
	\end{equation}
	
	By Lemma~\ref{lem:range}, $g_k\in\range\Hk$. Using \eqref{eq:v} and the fact that the positive eigenvalues of $\Hk$ lie in $[\mu,L_{g,1}]$, we have
	\begin{equation}
	\|v^k\|\leq\frac{1}{\mu}\|g_k\|.
	\label{eq:v-bound}
	\end{equation}
	To bound $r_k$, let $v^{t,k}$ be the $t$-th iterate of
	$\operatorname{GD}(\Hk,g_k;\gamma_k,T_1)$ and set
	$r^{t,k}:=g_k-\Hk v^{t,k}$. Then
	$r^{t+1,k}=(I-\gamma_k\Hk)r^{t,k}$.
	Since $r^{t,k}\in\range\Hk$, Assumption~\ref{asp:singularvalue} and
	$\gamma_k\leq1/L_{g,1}$ imply
	\begin{equation*}
	\|r^{t+1,k}\|
	\leq
	(1-\mu\gamma_k)\|r^{t,k}\|.
	\end{equation*}
	As $v^{0,k}=0$, after $T_1$ iterations,
	\begin{equation}
	\|r_k\|
	\leq
	(1-\mu\gamma_k)^{T_1}\|g_k\|
	\leq
	\frac{1}{4}\|g_k\|.
	\label{eq:rk}
	\end{equation}
	
	Combining \eqref{eq:y-descent0}, \eqref{eq:v-bound}, and \eqref{eq:rk},
	\begin{equation*}
	\|\nabla_2g(x^k,y^{k+1})\|
	\leq
	\frac{L_{g,2}}{2\mu^2}\|g_k\|^2+\frac{1}{4}\|g_k\|
	\leq\frac{1}{2}\|g_k\|,
	\end{equation*}
	where the last inequality uses
	$\|g_k\|\leq\mu^2/(2L_{g,2})$.
	Finally, the $L_{g,1}$-Lipschitz continuity of $\nabla_2g$ gives
	\begin{align*}
	\|\nabla_2g(x^{k+1},y^{k+1})\|
	\leq
	\|\nabla_2g(x^k,y^{k+1})\|
	+L_{g,1}\|x^{k+1}-x^k\|\leq
	\frac{1}{2}\|g_k\|
	+L_{g,1}\|x^{k+1}-x^k\|,
	\end{align*}
	which proves the result.
\end{proof} 

\begin{lemma}\label{lem:u}
	Suppose Assumptions~\ref{asp:singularvalue} and \ref{asp:smooth1} hold. In Newton-based GALET, if $0 \leq \eta_k\leq1/L_{g,1}$, then
	\begin{equation*}
	\big\|\Hk u^{k+1}
	-\Hk[\Hk]^\dagger\nabla_2f(x^k,y^k)\big\|
	\leq
	L_{f,0}(1-\mu\eta_k)^{T_2}.
	\end{equation*}
\end{lemma}

\begin{proof}
	Let $u^{t,k}$ be the $t$-th iterate of
	$\operatorname{GD}(\Hk,\nabla_2f(x^k,y^k);\eta_k,T_2)$ and define
	$e^{t,k}:=\Hk u^{t,k}-\Hk[\Hk]^\dagger\nabla_2f(x^k,y^k)$.
	The GD iteration gives
	$e^{t+1,k}=(I-\eta_k\Hk)e^{t,k}$.
	Since $e^{t,k}\in\range\Hk$, the nonzero eigenvalues of $\Hk$ lie in
	$[\mu,L_{g,1}]$, and $\eta_k\leq1/L_{g,1}$,
	\begin{equation*}
	\|e^{t+1,k}\|
	\leq
	(1-\mu\eta_k)\|e^{t,k}\|.
	\end{equation*}
	Moreover, $u^{0,k}=0$ and
	$P_k=\Hk[\Hk]^\dagger$ is the orthogonal projector onto $\range\Hk$, so
	\begin{equation*}
	\|e^{0,k}\|
	=
	\|P_k\nabla_2f(x^k,y^k)\|
	\leq
	\|\nabla_2f(x^k,y^k)\|
	\leq L_{f,0}.
	\end{equation*}
	Therefore,
	$\|e^{T_2,k}\|\leq L_{f,0}(1-\mu\eta_k)^{T_2}$.
	Since $u^{k+1}=u^{T_2,k}$, the result follows.
\end{proof} 


\begin{proof}[Proof of Theorem~\ref{thm:galet} ]
	We first establish by induction the local-region bound
	\begin{equation}
	\|\nabla_2g(x^k,y^k)\|
	\leq\frac{\mu^2}{2L_{g,2}},
	\qquad k=0,\ldots,K.
	\label{eq:local-region}
	\end{equation} 
	By Lemma~\ref{lem:range},
	$\range{\nabla_{21}^2g(x^k,y^k)}\subseteq\range\Hk$.
	Using
	$\|u^*(x^k,y^k)\|\leq L_{f,0}/\mu$,
	$\|[\Hk]^\dagger\|_{\operatorname{op}}\leq1/\mu$,
	and Lemma~\ref{lem:u}, we obtain
	\begin{equation*}
	\|d_x^k\|
	\leq C_{f,0}
	+\frac{L_{g,1}L_{f,0}}{\mu}
	\bigl(1+(1-\mu\eta_k)^{T_2}\bigr)
	\leq C_{f,0}+\frac{2L_{g,1}L_{f,0}}{\mu}
	=:D.
	\end{equation*}
	If \eqref{eq:local-region} holds at iteration $k$,
	Lemma~\ref{lem:y} and $x^{k+1}-x^k=-\alpha d_x^k$ give
	\begin{equation*}
	\|\nabla_2g(x^{k+1},y^{k+1})\|
	\leq\frac{\mu^2}{4L_{g,2}}+L_{g,1}\alpha D.
	\end{equation*}
	Thus, $\alpha\leq\mu^2/(4L_{g,1}L_{g,2}D)$
	preserves \eqref{eq:local-region}. Together with the
	initialization condition, induction proves
	\eqref{eq:local-region}. This requirement is satisfied by
	$\alpha=\Theta(\kappa^{-4})$ with a sufficiently small constant factor,
	under the stated parameter-dependence convention.
	
	Let $F_k:=F(x^k,y^k;u^*(x^k,y^k))$ and define
	the Lyapunov function
	\begin{equation*}
	V_k:=F_k-f^*
	+c\|\nabla_2g(x^k,y^k)\|^2
	+\frac{L_{f,0}^2}{c\mu^2},
	\end{equation*}
	where $c>0$ will be specified below.
	We first show that $V_k\geq0$.
	
	Let $\hat y^k$ be the Euclidean projection of $y^k$ onto
	$Y^*(x^k)$. By Assumption~\ref{asp:smooth1}(d) and
	Lemma~\ref{lem:range},
	\begin{equation*}
	f(x^k,y^k)
	\geq f(x^k,\hat y^k)
	-\frac{L_{f,0}}{\mu}\|\nabla_2g(x^k,y^k)\|
	\geq f^*-\frac{L_{f,0}}{\mu}\|\nabla_2g(x^k,y^k)\|.
	\end{equation*}
	Moreover, since $\|u^*(x^k,y^k)\|\leq L_{f,0}/\mu$,
	\begin{equation*}
	F_k
	\geq f(x^k,y^k)
	-\frac{L_{f,0}}{\mu}\|\nabla_2g(x^k,y^k)\|
	\geq f^*-\frac{2L_{f,0}}{\mu}
	\|\nabla_2g(x^k,y^k)\|.
	\end{equation*}
	Consequently, completing the square yields
	\begin{equation*}
	V_k
	\geq c\left(
	\|\nabla_2g(x^k,y^k)\|
	-\frac{L_{f,0}}{c\mu}
	\right)^2
	\geq0.
	\end{equation*}
	This lower bound holds globally, without requiring
	\eqref{eq:local-region}, which is needed for the
	Newton-type contraction analysis.
	
	Lemma~\ref{lem:y} gives
	\begin{equation}
	\|\nabla_2g(x^{k+1},y^{k+1})\|^2
	\leq
	\frac{1}{2}\|\nabla_2g(x^k,y^k)\|^2
	+2L_{g,1}^2\alpha^2\|d_x^k\|^2.
	\label{eq:descent-y}
	\end{equation}
	Using $\|v^k\|\leq\mu^{-1}\|\nabla_2g(x^k,y^k)\|$
	from \eqref{eq:v-bound}, Lemmas~\ref{lem:F} and
	\ref{lem:u} together with \eqref{eq:descent-y} yield
	\begin{align*}
	V_{k+1}-V_k
	\leq{}&
	\left(-\alpha+\alpha^2\left(
	\frac{1+L_F}{2}+\mu^2L_{u^*}^2
	+\frac{3L_{g,1}^2}{2\mu^2}
	+2cL_{g,1}^2\right)\right)\|d_x^k\|^2\\
	&+\left(-\frac{c}{2}+\frac{3+L_F}{2\mu^2}
	+L_{u^*}^2+\frac{3L_{g,1}^2}{2\mu^4}\right)
	\|\nabla_2g(x^k,y^k)\|^2
	+\frac{L_{g,1}^2L_{f,0}^2}{2\mu^2}
	(1-\mu\eta_k)^{2T_2}.
	\end{align*}
	The second coefficient is at most $-c/4$ if
	\begin{equation*}
	c\geq\frac{2(3+L_F)}{\mu^2}
	+4L_{u^*}^2+\frac{6L_{g,1}^2}{\mu^4},
	\end{equation*}
	while the first is at most $-\alpha/2$ if
	\begin{equation*}
	\alpha\leq
	\left(
	1+L_F+2\mu^2L_{u^*}^2
	+\frac{3L_{g,1}^2}{\mu^2}
	+4cL_{g,1}^2
	\right)^{-1}.
	\end{equation*}
	The bounds on $L_F$ and $L_{u^*}$, together with the
	parameter-dependence convention, allow us to choose
	$c=\Theta(L^2/\mu^4)$ with a sufficiently large constant
	and $\alpha=\Theta(\kappa^{-4})$ with a sufficiently
	small constant. These choices satisfy both descent
	conditions and preserve \eqref{eq:local-region}. 
	
	Since $\eta_k=\Theta(1/L_{g,1})$ and
	$L_{g,1}/\mu\leq\kappa$, there exists $c_\eta>0$,
	independent of $k$ and $\kappa$, such that
	$\mu\eta_k\geq c_\eta/\kappa$.
	Choosing $T_2\geq\kappa\log(\kappa^2K)/(2c_\eta)$ gives
	\begin{equation}
	(1-\mu\eta_k)^{2T_2}
	\leq e^{-2c_\eta T_2/\kappa}
	\leq\frac{1}{\kappa^2K}.
	\label{eq:u-linear-rate}
	\end{equation}
	Thus, $T_2=\Theta(\kappa\log(\kappa K))$ with a sufficiently
	large constant factor yields
	\begin{equation}
	V_{k+1}-V_k
	\leq
	-\frac{\alpha}{2}\|d_x^k\|^2
	-\frac{c}{4}\|\nabla_2g(x^k,y^k)\|^2
	+\frac{L_{f,0}^2}{2K},
	\label{eq:V-descent}
	\end{equation}
	where we used $L_{g,1}^2/\mu^2\leq\kappa^2$.
	
	It remains to bound $V_0$. By its definition,
	$\|u^*(x^0,y^0)\|\leq L_{f,0}/\mu$,
	and the initialization condition,
	\begin{align*}
	V_0
	&\leq f(x^0,y^0)-f^*
	+\frac{L_{f,0}}{\mu}\|\nabla_2g(x^0,y^0)\|
	+c\|\nabla_2g(x^0,y^0)\|^2
	+\frac{L_{f,0}^2}{c\mu^2}\\
	&\leq f(x^0,y^0)-f^*
	+\frac{L_{f,0}}{\kappa}
	+\frac{c\mu^4}{L^2}
	+\frac{L_{f,0}^2}{c\mu^2}.
	\end{align*}
	Since $c=\Theta(L^2/\mu^4)$ and $\kappa=L/\mu\geq1$,
	we have $V_0=O(1)$ with respect to $\kappa$,
	provided the initial objective gap and other problem constants
	are uniformly bounded.
	
	Summing \eqref{eq:V-descent} over $k=0,\ldots,K-1$
	and using $V_K\geq0$ gives
	\begin{align}
	\frac{1}{K}\sum_{k=0}^{K-1}\|d_x^k\|^2
	&=O\left(\frac{\kappa^4}{K}\right),
	\label{eq:Rx-rate}\\
	\frac{1}{K}\sum_{k=0}^{K-1}
	\|\nabla_2g(x^k,y^k)\|^2
	&=O\left(\frac{1}{cK}\right).
	\label{eq:gy-rate}
	\end{align}
	Since $\mathcal{R}_x(x^k,y^k,u^{k+1})=\|d_x^k\|^2$,
	\eqref{eq:Rx-rate} gives the claimed $\mathcal{R}_x$ bound.
	By the PL inequality following Lemma~\ref{lem:range},
	$\mathcal{R}_y(x^k,y^k)
	\leq\mu^{-1}\|\nabla_2g(x^k,y^k)\|^2$.
	Thus, \eqref{eq:gy-rate} and $c=\Theta(L^2/\mu^4)$ yield
	\begin{equation*}
	\frac{1}{K}\sum_{k=0}^{K-1}\mathcal{R}_y(x^k,y^k)
	=O\left(\frac{\mu}{\kappa^2K}\right)
	=O\left(\frac{1}{\kappa^2K}\right),
	\end{equation*}
	where the last step uses $0<\mu\leq1$.
	
	For the multiplier residual, $\Hk(I-P_k)=0$ and
	Lemma~\ref{lem:u} give
	\begin{align*}
	\mathcal{R}_u(x^k,y^k,u^{k+1})
	&=\|\Hk(P_k\nabla_2f(x^k,y^k)-\Hk u^{k+1})\|^2\\
	&\leq L_{g,1}^2L_{f,0}^2(1-\mu\eta_k)^{2T_2}.
	\end{align*}
	By \eqref{eq:u-linear-rate} and $L_{g,1}/\mu\leq\kappa$,
	we obtain
	\begin{equation*}
	\frac{1}{K}\sum_{k=0}^{K-1}
	\mathcal{R}_u(x^k,y^k,u^{k+1})
	=O(K^{-1}).
	\end{equation*}
	
	Combining the three bounds, there exists
	$k^*\in\{0,\ldots,K-1\}$ such that all three residuals at
	$(x^{k^*},y^{k^*},u^{k^*+1})$ are $O(\kappa^4/K)$.
	Thus, $K=O(\kappa^4/\epsilon)$ outer iterations suffice for
	$\epsilon$-stationarity. Each outer iteration requires $O(1)$ gradient/JVP computations
	and $T_1+T_2=O(\kappa\log(\kappa K))$ Hessian-vector products.
	With $K=O(\kappa^4/\epsilon)$, the total complexities are
	$O(\kappa^4/\epsilon)$ gradient/JVP computations and
	$O(\kappa^5\epsilon^{-1}\log(\kappa/\epsilon))$
	Hessian-vector products.
\end{proof}


\subsection{Convergence Analysis for Newton-Based BAMM}

For convenience, define
\begin{align}
L_{\psi,1}
&:=L_{g,1}+L_{f,1}\sup_{k\geq0}\sigma_k,
\qquad
L_{\psi,2}
:=L_{g,2}+L_{f,2}\sup_{k\geq0}\sigma_k, \nonumber\\
C_{\Phi 1}
&:=L_{\psi,1}^2L_{\psi,2}
+2L_{\psi,1}L_{\psi,2}\mu_f\sup_{k\geq0}\sigma_k
+L_{\psi,2}\mu_f^2\sup_{k\geq0}\sigma_k^2,
\label{eq:constant-phi}\\
C_{\Phi 2}
&:=L_{\psi,1}^2L_{f,1}
+2L_{\psi,1}L_{f,1}\mu_f\sup_{k\geq0}\sigma_k
+L_{f,1}\mu_f^2\sup_{k\geq0}\sigma_k^2.
\nonumber
\end{align}
The following lemma gives one-step tracking bounds for the LL and multiplier variables.

\begin{lemma}\label{lem:bamm-tracking}
	Suppose Assumptions~\ref{asp:fsc} and \ref{asp:smooth2} hold. In Newton-based BAMM, let
	$0 \leq \gamma_k\leq1/L_{\psi,1}$ and
	$T_k\geq\log(1/4)/\log(1-\mu_f\sigma_k\gamma_k)$.
	Then
	\begin{align}
	\|y^{k+1}-y_{\sigma_{k+1}}^*(x^{k+1})\|
	\leq{}&
	\frac{1}{4}\|y^k-y_{\sigma_k}^*(x^k)\|
	+\frac{L_{\psi,2}}{2\mu_f\sigma_k}
	\|y^k-y_{\sigma_k}^*(x^k)\|^2 \label{eq:y-descent1}\\
	&+\frac{L_{\psi,1}}{\mu_f\sigma_k}\|x^{k+1}-x^k\|
	+\frac{2\|\nabla_2f(x^{k+1},y_{\sigma_{k+1}}^*(x^{k+1}))\|}{\mu_f}
	\frac{\sigma_k-\sigma_{k+1}}{\sigma_k}, \nonumber
	\end{align}
	and
	\begin{align}
	&\|u^{k+1}-u_{\sigma_{k+1}}^*(x^{k+1})\| \nonumber \\
	\leq{}&
	\frac{1}{4}\|u^k-u_{\sigma_k}^*(x^k)\|
	+\frac{5\bigl(L_{f,1}+L_{\psi,2}\|u_{\sigma_k}^*(x^k)\|\bigr)}
	{4\mu_f\sigma_k}
	\|y^k-y_{\sigma_k}^*(x^k)\| \nonumber\\
	&+\left(1+\frac{L_{\psi,1}}{\mu_f\sigma_k}\right)
	\frac{L_{f,1}+L_{\psi,2}\|u_{\sigma_k}^*(x^k)\|}
	{\mu_f\sigma_k}
	\|x^{k+1}-x^k\| \label{eq:u-descent1}\\
	&+\left(
	\frac{(2L_{f,2}+2L_{g,2})\|u_{\sigma_{k+1}}^*(x^{k+1})\|}{\mu_f^2}
	+\frac{4L_{f,1}+2L_{g,1}}{\mu_f^2}
	\right)
	\|\nabla_2f(x^{k+1},y_{\sigma_{k+1}}^*(x^{k+1}))\|
	\frac{\sigma_k-\sigma_{k+1}}{\sigma_k^2}.
	\nonumber
	\end{align}
\end{lemma}

\begin{proof}
	Let $y_k^*:=y_{\sigma_k}^*(x^k)$ and
	$H_k^\psi:=\nabla_{22}^2\psi_{\sigma_k}(x^k,y^k)$.
	Let $v^{t,k}$ be the $t$-th iterate of
	$\operatorname{GD}(H_k^\psi,\nabla_2\psi_{\sigma_k}(x^k,y^k);
	\gamma_k,T_k)$.
	As in \cite[Lemma~E.3]{dong2025efficient},
	\begin{equation*}
	\|y^k-v^{t+1,k}-y_k^*\|
	\leq
	(1-\mu_f\sigma_k\gamma_k)
	\|y^k-v^{t,k}-y_k^*\|
	+\frac{L_{\psi,2}\gamma_k}{2}\|y^k-y_k^*\|^2.
	\end{equation*}
	After $T_k$ iterations,
	\begin{align}
	\|y^{k+1}-y_k^*\|
	&\leq
	(1-\mu_f\sigma_k\gamma_k)^{T_k}
	\|y^k-y_k^*\|
	+\frac{L_{\psi,2}}{2\mu_f\sigma_k}
	\|y^k-y_k^*\|^2 \nonumber\\
	&\leq
	\frac{1}{4}\|y^k-y_k^*\|
	+\frac{L_{\psi,2}}{2\mu_f\sigma_k}
	\|y^k-y_k^*\|^2.
	\label{eq:y-descent}
	\end{align}
	Moreover, \cite[Lemma~D.8]{liu2023averaged} gives
	\begin{align}
	\|y_k^*-y_{\sigma_{k+1}}^*(x^{k+1})\|
	\leq{}
	\frac{L_{\psi,1}}{\mu_f\sigma_k}
	\|x^{k+1}-x^k\| +
	\frac{2\|\nabla_2f(x^{k+1},y_{\sigma_{k+1}}^*(x^{k+1}))\|}{\mu_f}
	\frac{\sigma_k-\sigma_{k+1}}{\sigma_k}.
	\label{eq:ysssk-error}
	\end{align}
	Combining \eqref{eq:y-descent} and \eqref{eq:ysssk-error} proves \eqref{eq:y-descent1}.
	
	For the multiplier update, define
	$u_{\sigma_k}^*(x^k,y^k)
	:=(H_k^\psi)^{-1}\nabla_2f(x^k,y^k)$.
	Let $w^{t,k}$ be the $t$-th iterate of $\operatorname{GD}\bigl(
	H_k^\psi,
	H_k^\psi u^k-\nabla_2f(x^k,y^k);
	\gamma_k,T_k\bigr).$
	As in \cite[Lemma~E.3]{dong2025efficient},
	\begin{equation*}
	\|u^k-w^{t+1,k}-u_{\sigma_k}^*(x^k,y^k)\|
	\leq
	(1-\mu_f\sigma_k\gamma_k)
	\|u^k-w^{t,k}-u_{\sigma_k}^*(x^k,y^k)\|.
	\end{equation*}
	Since $w^{0,k}=0$, after $T_k$ iterations,
	\begin{align*}
	\|u^{k+1}-u_{\sigma_k}^*(x^k,y^k)\|
	\leq{}&
	\frac{1}{4}\|u^k-u_{\sigma_k}^*(x^k)\|
	+\frac{1}{4}
	\|u_{\sigma_k}^*(x^k,y^k)-u_{\sigma_k}^*(x^k)\|.
	\end{align*}
	Consequently,
	\begin{align}
	\|u^{k+1}-u_{\sigma_{k+1}}^*(x^{k+1})\|
	\leq{}&
	\frac{1}{4}\|u^k-u_{\sigma_k}^*(x^k)\|
	+\frac{5}{4}
	\|u_{\sigma_k}^*(x^k,y^k)-u_{\sigma_k}^*(x^k)\| \nonumber\\
	&+
	\|u_{\sigma_k}^*(x^k)-u_{\sigma_{k+1}}^*(x^{k+1})\|.
	\label{eq:u-descent2}
	\end{align}
	Using the $(\mu_f\sigma_k)$-strong convexity of
	$\psi_{\sigma_k}(x^k,\cdot)$ and an argument similar to that for
	\citep[Eq.~(21)]{dong2025efficient}, we obtain
	\begin{equation}
	\|u_{\sigma_k}^*(x^k,y^k)-u_{\sigma_k}^*(x^k)\|
	\leq
	\frac{L_{f,1}+L_{\psi,2}\|u_{\sigma_k}^*(x^k)\|}
	{\mu_f\sigma_k}
	\|y^k-y_k^*\|.
	\label{eq:u-star}
	\end{equation} 
	Finally, \cite[Lemma~D.8]{liu2023averaged} gives
	\begin{align}
	& \|u_{\sigma_{k+1}}^*(x^{k+1})-u_{\sigma_k}^*(x^k)\| \nonumber \\
	\leq{}&
	\left(1+\frac{L_{\psi,1}}{\mu_f\sigma_k}\right)
	\frac{L_{f,1}+L_{\psi,2}\|u_{\sigma_k}^*(x^k)\|}
	{\mu_f\sigma_k}
	\|x^{k+1}-x^k\| 
	\label{eq:uss}\\
	&+
	\left(
	\frac{(2L_{f,2}+2L_{g,2})\|u_{\sigma_{k+1}}^*(x^{k+1})\|}{\mu_f^2}
	+\frac{4L_{f,1}+2L_{g,1}}{\mu_f^2}
	\right)
	\|\nabla_2f(x^{k+1},y_{\sigma_{k+1}}^*(x^{k+1}))\|
	\frac{\sigma_k-\sigma_{k+1}}{\sigma_k^2}.\nonumber
	\end{align}
	Combining \eqref{eq:u-descent2}, \eqref{eq:u-star}, and \eqref{eq:uss} proves \eqref{eq:u-descent1}.
\end{proof} 


\begin{proof}[Proof of Theorem~\ref{thm:bamm}]
	Let $\Phi_\sigma(x):=f(x,y_\sigma^*(x))$,
	$e_y^k:=y^k-y_{\sigma_k}^*(x^k)$, and
	$e_u^k:=u^k-u_{\sigma_k}^*(x^k)$.
	Define
	\begin{equation*}
	V_k:=
	a_k\bigl(\Phi_{\sigma_k}(x^k)-f_0\bigr)
	+b_k\|e_y^k\|^2+c_k\|e_u^k\|^2,
	\end{equation*}
	where $a_k,b_k,c_k>0$ are non-increasing.
	Since $\Phi_{\sigma_k}(x^k)\geq f_0$, we have $V_k\geq0$.
	By \cite[Lemma~C.1]{liu2023averaged},
	\begin{align}
	&\Phi_{\sigma_{k+1}}(x^{k+1})-\Phi_{\sigma_k}(x^k)
	\leq
	-\frac{\alpha_k}{2}\|\nabla\Phi_{\sigma_k}(x^k)\|^2
	\nonumber\\
	&\quad-\frac12
	\left(
	\frac1{\alpha_k}
	-\frac{C_{\Phi 1}\|u_{\sigma_k}^*(x^k)\|+C_{\Phi 2}}
	{\mu_f^2\sigma_k^2}
	\right)\|x^{k+1}-x^k\|^2
	+\alpha_kL_{\psi,1}^2\|e_u^k\|^2
	\nonumber\\
	&\quad+\alpha_k
	\bigl(L_{\psi,2}\|u_{\sigma_k}^*(x^k)\|+L_{f,1}\bigr)^2
	\|e_y^k\|^2
	+\frac{2\|\nabla_2f(x^{k+1},
		y_{\sigma_{k+1}}^*(x^{k+1}))\|^2}{\mu_f}
	\frac{\sigma_k-\sigma_{k+1}}{\sigma_k}.
	\label{eq:phis-descent}
	\end{align}
	
	Let $U_k:=\|u_{\sigma_k}^*(x^k)\|$ and
	$G_{k+1}:=\|\nabla_2f(x^{k+1},
	y_{\sigma_{k+1}}^*(x^{k+1}))\|$.
	Whenever $\|e_y^k\|\leq\mu_f\sigma_k/(2L_{\psi,2})$,
	combining \eqref{eq:phis-descent} with
	\eqref{eq:y-descent1}--\eqref{eq:u-descent1}
	and the monotonicity of the Lyapunov weights gives
	\begin{align}
	&V_{k+1}-V_k
	\leq
	-\frac{a_k\alpha_k}{2}
	\|\nabla\Phi_{\sigma_k}(x^k)\|^2
	\nonumber\\
	&-\left(
	\frac{a_k}{2\alpha_k}
	-\frac{a_k(C_{\Phi 1}U_k+C_{\Phi 2})+8b_kL_{\psi,1}^2}
	{2\mu_f^2\sigma_k^2}
	-\frac{6c_k(L_{f,1}+L_{\psi,2}U_k)^2}
	{\mu_f^2\sigma_k^2}
	\left(1+\frac{L_{\psi,1}}{\mu_f\sigma_k}\right)^2
	\right)\|x^{k+1}-x^k\|^2
	\nonumber\\
	&-\left(
	\frac{b_k}{2}
	-a_k\alpha_k(L_{\psi,2}U_k+L_{f,1})^2
	-\frac{75c_k(L_{f,1}+L_{\psi,2}U_k)^2}
	{8\mu_f^2\sigma_k^2}
	\right)\|e_y^k\|^2
	\nonumber\\
	&-\left(
	\frac{7c_k}{8}-a_k\alpha_kL_{\psi,1}^2
	\right)\|e_u^k\|^2
	\nonumber\\
	&+\frac{2a_kG_{k+1}^2}{\mu_f}
	\frac{\sigma_k-\sigma_{k+1}}{\sigma_k}
	+\frac{16b_kG_{k+1}^2}{\mu_f^2}
	\left(\frac{\sigma_k-\sigma_{k+1}}{\sigma_k}\right)^2
	\nonumber\\
	&+\frac{6c_kG_{k+1}^2}{\mu_f^4}
	\bigl((2L_{f,2}+2L_{g,2})U_{k+1}
	+4L_{f,1}+2L_{g,1}\bigr)^2
	\left(\frac{\sigma_k-\sigma_{k+1}}{\sigma_k^2}\right)^2.
	\label{eq:bamm-V-descent}
	\end{align}
	
	Define
	\begin{equation*}
	\Pi_k:=
	\|\nabla\Phi_{\sigma_k}(x^k)\|^2
	+\|x^{k+1}-x^k\|^2
	+\|e_y^k\|^2+\|e_u^k\|^2.
	\end{equation*}
	We next verify the local tracking conditions and Lyapunov
	descent for each boundedness setting, with $k_0$ chosen
	independently of $K$.
	
	\paragraph{Case 1}
	Suppose $\{x^k\}$, $\{y^k\}$, and $\{u^k\}$ are bounded.
	By \cite[Lemma~D.5]{liu2023averaged},
	\begin{equation*}
	\|\nabla_2f(x^k,y_{\sigma_k}^*(x^k))\|
	=O(\sigma_k^{-1}),
	\qquad
	\|u_{\sigma_k}^*(x^k)\|=O(\sigma_k^{-2}).
	\end{equation*}
	Choose
	\begin{equation*}
	\sigma_k=\frac{\bar{\sigma}}{(k+k_0)^p},\quad
	\alpha_k=\frac{\bar{\alpha}}{(k+k_0)^{4p+2\tau}},\quad
	a_k=\frac{1}{(k+k_0)^{2p+\tau}},\quad
	b_k=1,\quad
	c_k=\frac{1}{(k+k_0)^{6p+\tau}}.
	\end{equation*}
	
	We first verify the local condition in \eqref{eq:bamm-V-descent}.
	By Assumption~\ref{asp:smooth2},
	$\|d_x^k\|\leq C_{f,0}+L_{\psi,1}(\|e_u^k\|+U_k)$.
	Since $U_k=O((k+k_0)^{2p})$,
	\begin{equation*}
	\|x^{k+1}-x^k\|
	=O((k+k_0)^{-4p-2\tau})\|e_u^k\|
	+O((k+k_0)^{-2p-2\tau}).
	\end{equation*}
	Moreover,
	$(\sigma_k-\sigma_{k+1})/\sigma_k
	=O((k+k_0)^{-1})$ and
	$(\sigma_k-\sigma_{k+1})/\sigma_k^2
	=O((k+k_0)^{-1+p})$.
	Whenever $\|e_y^k\|\leq\mu_f\sigma_k/(2L_{\psi,2})$,
	\eqref{eq:y-descent1}--\eqref{eq:u-descent1} give
	\begin{align}
	\|e_y^{k+1}\|
	&\leq \frac12\|e_y^k\|
	+O((k+k_0)^{-3p-2\tau})\|e_u^k\|
	+O((k+k_0)^{-p-2\tau})
	+O((k+k_0)^{-1+p}),
	\label{eq:case1-y}\\
	\|e_u^{k+1}\|
	&\leq \frac14\|e_u^k\|
	+O((k+k_0)^{3p})\|e_y^k\|
	+O((k+k_0)^{-2\tau})\|e_u^k\|
	+O((k+k_0)^{2p-2\tau})
	\nonumber\\
	&\quad+O((k+k_0)^{-1+4p}).
	\label{eq:case1-u}
	\end{align}
	
	We prove by induction that
	\begin{equation}
	\|e_y^k\|\leq\frac{\mu_f\sigma_k}{2L_{\psi,2}},
	\qquad
	\|e_u^k\|\leq(k+k_0)^{2p+\tau}.
	\label{eq:case1-induction}
	\end{equation}
	Under \eqref{eq:case1-induction}, the perturbations in
	\eqref{eq:case1-y} are
	$O((k+k_0)^{-p-2\tau})$,
	$O((k+k_0)^{-p-\tau})$, and
	$O((k+k_0)^{-1+p})$,
	all $o(\sigma_k)$ for $p<1/2$.
	Thus, for sufficiently large $k_0$,
	$\|e_y^{k+1}\|
	\leq\mu_f\sigma_{k+1}/(2L_{\psi,2})$.
	
	Similarly, \eqref{eq:case1-u} gives
	\begin{equation*}
	\|e_u^{k+1}\|
	\leq
	\frac14(k+k_0)^{2p+\tau}
	+O((k+k_0)^{2p})
	+O((k+k_0)^{2p-2\tau})
	+O((k+k_0)^{2p-\tau})
	+O((k+k_0)^{-1+4p}),
	\end{equation*}
	which is at most $(k+k_0+1)^{2p+\tau}$
	for sufficiently large $k_0$ and $p<1/2$.
	The initial LL bound holds by assumption, while
	$\|e_u^0\|\leq\|u^0\|+O(k_0^{2p})
	\leq k_0^{2p+\tau}$ for sufficiently large $k_0$.
	Hence \eqref{eq:case1-induction} holds for all $k$.
	
	Substituting these bounds into \eqref{eq:bamm-V-descent}
	and choosing $k_0$ sufficiently large gives
	\begin{align}
	V_{k+1}-V_k
	\leq{}&
	-\frac{a_k\alpha_k}{2}
	\|\nabla\Phi_{\sigma_k}(x^k)\|^2
	-\frac{a_k}{4\alpha_k}\|x^{k+1}-x^k\|^2
	-\frac14\|e_y^k\|^2
	-\frac{7c_k}{16}\|e_u^k\|^2
	\nonumber\\
	&+O((k+k_0)^{-1-\tau})
	+O((k+k_0)^{-2+2p})
	+O((k+k_0)^{-2-\tau+2p}).
	\label{eq:case1-V}
	\end{align}
	For $p<1/2$ and $\tau>0$, the remainder terms are summable.
	Since $a_k\alpha_k=\Theta((k+k_0)^{-6p-3\tau})$
	is, up to constants, the smallest descent weight,
	telescoping \eqref{eq:case1-V} and using $V_K\geq0$ give
	\begin{equation}
	\frac1K\sum_{k=0}^{K-1}\Pi_k
	=O\left(\frac{1}{Ka_K\alpha_K}\right)
	=O\left(K^{-(1-6p-3\tau)}\right).
	\label{eq:case1-Pi}
	\end{equation}
	By the boundedness of $\{u^k\}$ and
	\cite[Eq.~(35)]{liu2023averaged},
	\begin{equation*}
	\min_{0\leq k\leq K-1}
	\operatorname{KKT}(x^k,y^k,u^k)
	=O\left(K^{-(1-6p-3\tau)}+K^{-2p}\right).
	\end{equation*}
	Balancing the exponents gives $p=(1-3\tau)/8$,
	and hence an $\epsilon$-stationary iterate within
	$K=O(\epsilon^{-4/(1-3\tau)})$ outer iterations.
	
	Each outer iteration requires $O(1)$ gradient/JVP computations
	and $T_k=\Theta((k+k_0)^p)$ Hessian-vector products.
	Since $k_0$ is independent of $K$,
	$\sum_{k=0}^{K-1}T_k=O(K^{1+p})$.
	With $p=(1-3\tau)/8$, the total complexities are
	$O(\epsilon^{-4/(1-3\tau)})$ gradient/JVP computations and
	$O(\epsilon^{-(9-3\tau)/(2-6\tau)})$ Hessian-vector products. 
	
	\paragraph{Case 2}
	Suppose $\{\nabla_2f(x^k,y_{\sigma_k}^*(x^k))\}$ and
	$\{u^k\}$ are bounded. By the $(\mu_f\sigma_k)$-strong
	convexity of $\psi_{\sigma_k}(x^k,\cdot)$,
	\begin{equation*}
	\|u_{\sigma_k}^*(x^k)\|
	\leq
	\frac{\|\nabla_2f(x^k,y_{\sigma_k}^*(x^k))\|}
	{\mu_f\sigma_k}
	=O(\sigma_k^{-1}).
	\end{equation*}
	Choose
	\begin{equation*}
	\sigma_k=\frac{\bar{\sigma}}{(k+k_0)^p},\quad
	\alpha_k=\frac{\bar{\alpha}}{(k+k_0)^{3p+2\tau}},\quad
	a_k=\frac{1}{(k+k_0)^{p+\tau}},\quad
	b_k=1,\quad
	c_k=\frac{1}{(k+k_0)^{4p+\tau}}.
	\end{equation*}
	
	By Assumption~\ref{asp:smooth2} and the choice of $\alpha_k$,
	\begin{equation*}
	\|x^{k+1}-x^k\|
	=O((k+k_0)^{-3p-2\tau})\|e_u^k\|
	+O((k+k_0)^{-2p-2\tau}).
	\end{equation*}
	Whenever $\|e_y^k\|\leq\mu_f\sigma_k/(2L_{\psi,2})$,
	\eqref{eq:y-descent1}--\eqref{eq:u-descent1} give
	\begin{align}
	\|e_y^{k+1}\|
	&\leq
	\frac12\|e_y^k\|
	+O((k+k_0)^{-2p-2\tau})\|e_u^k\|
	+O((k+k_0)^{-p-2\tau})
	+O((k+k_0)^{-1}),
	\label{eq:case2-y}\\
	\|e_u^{k+1}\|
	&\leq
	\frac14\|e_u^k\|
	+O((k+k_0)^{2p})\|e_y^k\|
	+O((k+k_0)^{-2\tau})\|e_u^k\|
	\nonumber\\
	&\quad+O((k+k_0)^{p-2\tau})
	+O((k+k_0)^{-1+2p}).
	\label{eq:case2-u}
	\end{align}
	
	We prove by induction that
	\begin{equation}
	\|e_y^k\|\leq\frac{\mu_f\sigma_k}{2L_{\psi,2}},
	\qquad
	\|e_u^k\|\leq(k+k_0)^{p+\tau}.
	\label{eq:case2-induction}
	\end{equation}
	Under \eqref{eq:case2-induction}, the perturbations in
	\eqref{eq:case2-y} are $o(\sigma_k)$ for $0<p<1$,
	so the LL bound holds at iteration $k+1$ for sufficiently
	large $k_0$. Similarly,
	\begin{equation*}
	\|e_u^{k+1}\|
	\leq
	\frac14(k+k_0)^{p+\tau}
	+O((k+k_0)^p)
	+O((k+k_0)^{p-\tau})
	+O((k+k_0)^{-1+2p})
	\leq(k+k_0+1)^{p+\tau}.
	\end{equation*}
	The initial LL bound holds by assumption, while
	$\|e_u^0\|\leq\|u^0\|+O(k_0^p)
	\leq k_0^{p+\tau}$ for sufficiently large $k_0$.
	Hence \eqref{eq:case2-induction} holds for all $k$.
	
	Substituting these bounds into \eqref{eq:bamm-V-descent}
	and choosing $k_0$ sufficiently large gives
	\begin{align}
	V_{k+1}-V_k
	\leq{}&
	-\frac{a_k\alpha_k}{2}
	\|\nabla\Phi_{\sigma_k}(x^k)\|^2
	-\frac{a_k}{4\alpha_k}\|x^{k+1}-x^k\|^2
	-\frac14\|e_y^k\|^2
	-\frac{7c_k}{16}\|e_u^k\|^2
	\nonumber\\
	&+O((k+k_0)^{-1-p-\tau})
	+O((k+k_0)^{-2})
	+O((k+k_0)^{-2-\tau}).
	\label{eq:case2-V}
	\end{align}
	The remainder terms are summable. Since
	$a_k\alpha_k=\Theta((k+k_0)^{-4p-3\tau})$
	is, up to constants, the smallest descent weight,
	telescoping \eqref{eq:case2-V} and using $V_K\geq0$ give
	\begin{equation}
	\frac1K\sum_{k=0}^{K-1}\Pi_k
	=O\left(\frac{1}{Ka_K\alpha_K}\right)
	=O\left(K^{-(1-4p-3\tau)}\right).
	\label{eq:case2-Pi}
	\end{equation}
	By the boundedness of $\{u^k\}$ and
	\cite[Eq.~(35)]{liu2023averaged},
	\begin{equation*}
	\min_{0\leq k\leq K-1}
	\operatorname{KKT}(x^k,y^k,u^k)
	=O\left(K^{-(1-4p-3\tau)}+K^{-2p}\right).
	\end{equation*}
	Balancing the exponents gives $p=(1-3\tau)/6$,
	and hence an $\epsilon$-stationary iterate within
	$K=O(\epsilon^{-3/(1-3\tau)})$ outer iterations.
	
	Each outer iteration requires $O(1)$ gradient/JVP computations
	and $T_k=\Theta((k+k_0)^p)$ Hessian-vector products.
	Since $k_0$ is independent of $K$,
	$\sum_{k=0}^{K-1}T_k=O(K^{1+p})$.
	With $p=(1-3\tau)/6$, the total complexities are
	$O(\epsilon^{-3/(1-3\tau)})$ gradient/JVP computations and
	$O(\epsilon^{-(7-3\tau)/(2-6\tau)})$ Hessian-vector products. 
	
	\paragraph{Case 3}
	Suppose $\{u_{\sigma_k}^*(x^k)\}$ and $\{u^k\}$ are bounded.
	Since
	\begin{equation*}
	\nabla_2f(x^k,y_{\sigma_k}^*(x^k))
	=
	\nabla_{22}^2\psi_{\sigma_k}(x^k,y_{\sigma_k}^*(x^k))
	u_{\sigma_k}^*(x^k),
	\end{equation*}
	we have
	$\|\nabla_2f(x^k,y_{\sigma_k}^*(x^k))\|
	\leq L_{\psi,1}\|u_{\sigma_k}^*(x^k)\|=O(1)$.
	Choose
	\begin{equation*}
	\sigma_k=\frac{\bar{\sigma}}{(k+k_0)^p},\quad
	\alpha_k=\frac{\bar{\alpha}}{(k+k_0)^{2p+2\tau}},\quad
	a_k=\frac{1}{(k+k_0)^\tau},\quad
	b_k=1,\quad
	c_k=\frac{1}{(k+k_0)^{2p+\tau}}.
	\end{equation*}
	
	By Assumption~\ref{asp:smooth2} and the choice of $\alpha_k$,
	\begin{equation*}
	\|x^{k+1}-x^k\|
	=O((k+k_0)^{-2p-2\tau})(1+\|e_u^k\|).
	\end{equation*}
	Whenever $\|e_y^k\|\leq\mu_f\sigma_k/(2L_{\psi,2})$,
	\eqref{eq:y-descent1}--\eqref{eq:u-descent1} give
	\begin{align}
	\|e_y^{k+1}\|
	&\leq
	\frac12\|e_y^k\|
	+O((k+k_0)^{-p-2\tau})(1+\|e_u^k\|)
	+O((k+k_0)^{-1}),
	\label{eq:case3-y}\\
	\|e_u^{k+1}\|
	&\leq
	\frac14\|e_u^k\|
	+O((k+k_0)^p)\|e_y^k\|
	+O((k+k_0)^{-2\tau})(1+\|e_u^k\|)
	+O((k+k_0)^{-1+p}).
	\label{eq:case3-u}
	\end{align}
	
	We prove by induction that
	\begin{equation}
	\|e_y^k\|\leq\frac{\mu_f\sigma_k}{2L_{\psi,2}},
	\qquad
	\|e_u^k\|\leq(k+k_0)^\tau.
	\label{eq:case3-induction}
	\end{equation}
	Under \eqref{eq:case3-induction}, the perturbations in
	\eqref{eq:case3-y} are $o(\sigma_k)$ for $0<p<1$,
	so the LL bound holds at iteration $k+1$ for sufficiently
	large $k_0$. Similarly,
	\begin{equation*}
	\|e_u^{k+1}\|
	\leq
	\frac14(k+k_0)^\tau
	+O(1)+O((k+k_0)^{-\tau})
	+O((k+k_0)^{-1+p})
	\leq(k+k_0+1)^\tau.
	\end{equation*}
	The initial LL bound holds by assumption, while
	$\|e_u^0\|=O(1)\leq k_0^\tau$
	for sufficiently large $k_0$.
	Hence \eqref{eq:case3-induction} holds for all $k$.
	
	Substituting these bounds into \eqref{eq:bamm-V-descent}
	and choosing $k_0$ sufficiently large gives
	\begin{align}
	V_{k+1}-V_k
	\leq{}&
	-\frac{a_k\alpha_k}{2}
	\|\nabla\Phi_{\sigma_k}(x^k)\|^2
	-\frac{a_k}{4\alpha_k}\|x^{k+1}-x^k\|^2
	-\frac14\|e_y^k\|^2
	-\frac{7c_k}{16}\|e_u^k\|^2
	\nonumber\\
	&+O((k+k_0)^{-1-\tau})
	+O((k+k_0)^{-2})
	+O((k+k_0)^{-2-\tau}).
	\label{eq:case3-V}
	\end{align}
	The remainder terms are summable.
	Since $a_k\alpha_k=\Theta((k+k_0)^{-2p-3\tau})$
	is, up to constants, the smallest descent weight,
	telescoping \eqref{eq:case3-V} and using $V_K\geq0$ give
	\begin{equation}
	\frac1K\sum_{k=0}^{K-1}\Pi_k
	=O\left(\frac{1}{Ka_K\alpha_K}\right)
	=O\left(K^{-(1-2p-3\tau)}\right).
	\label{eq:case3-Pi}
	\end{equation}
	By the boundedness of $\{u^k\}$ and
	\cite[Eq.~(35)]{liu2023averaged},
	\begin{equation*}
	\min_{0\leq k\leq K-1}
	\operatorname{KKT}(x^k,y^k,u^k)
	=O\left(K^{-(1-2p-3\tau)}+K^{-2p}\right).
	\end{equation*}
	Balancing the exponents gives $p=(1-3\tau)/4$,
	and hence an $\epsilon$-stationary iterate within
	$K=O(\epsilon^{-2/(1-3\tau)})$ outer iterations.
	
	Each outer iteration requires $O(1)$ gradient/JVP computations
	and $T_k=\Theta((k+k_0)^p)$ Hessian-vector products.
	Since $k_0$ is independent of $K$,
	$\sum_{k=0}^{K-1}T_k=O(K^{1+p})$.
	With $p=(1-3\tau)/4$, the total complexities are
	$O(\epsilon^{-2/(1-3\tau)})$ gradient/JVP computations and
	$O(\epsilon^{-(5-3\tau)/(2-6\tau)})$ Hessian-vector products.
\end{proof} 


\subsection{Derivation of the GALET Baseline Complexity Bounds}
\label{app:galet-complexity}

We derive condition-number-dependent complexity bounds for GALET from the analysis in Appendix~F of \citet{xiao2023generalized}.

Following \citet{xiao2023generalized}, we take $N=\Theta(1)$. Under our parameter-dependence convention, the constants $L_w$, $L_F$, and $L_g$ in their analysis satisfy $L_w=O(L/\mu^2)$, $L_F=O(L/\mu)$, and $L_g=O(L^2/\mu)$. The parameter restrictions following Eq.~(42) then permit $c=\Theta(L^3/\mu^4)$ and $\beta=\Theta(1/L)$. In particular, the UL stepsize restriction $\alpha\leq\mu_g^2\beta/(8cL_{g,1}^2)$, with $\mu_g=\mu/2$, permits
\begin{equation*}
\alpha=\Theta\left(\frac{\mu^2\beta}{cL^2}\right)
=\Theta(\kappa^{-6}).
\end{equation*}

GALET uses the Lyapunov function $V^k=F(x^k,y^k;u^*(x^k,y^k)) +c\mathcal{R}_y(x^k,y^k)$.
We first establish a uniform lower bound. By the LL error bound and $\|u^*(x,y)\|\leq L_{f,0}/\mu$, as in the proof of Theorem~\ref{thm:galet},
\begin{equation*}
F(x,y;u^*(x,y))
\geq f^*-\frac{2L_{f,0}}{\mu}\|\nabla_2g(x,y)\|.
\end{equation*}
By smoothness of $g(x,\cdot)$, $\mathcal{R}_y\geq\|\nabla_2g\|^2/(2L_{g,1})$, 
completing the square yields 
\begin{equation*}
V^k
\geq f^*-\frac{2L_{f,0}^2L_{g,1}}{c\mu^2}
=f^*-O(\kappa^{-2}).
\end{equation*}
Thus, we may choose a lower bound $\underline V=f^*-C$, where $C=O(1)$ with respect to $\kappa$.

Summing \cite[Eq.~(42)]{xiao2023generalized} and using $V^K\geq\underline V$ gives 
\begin{equation*}
\frac{1}{K}\sum_{k=0}^{K-1}\|d_x^k\|^2
\leq
\frac{16(V^0-\underline V)}{\alpha K}
+\frac{8L_{g,1}^2}{K}\sum_{k=0}^{K-1}b_k^2,
\end{equation*}
where $b_k$ measures the multiplier approximation error. Following \citet{xiao2023generalized}, this error is controlled by choosing $T=\Theta(\kappa^2\log(\kappa K))$ inner iterations. The first term therefore depends on the initial Lyapunov gap $V^0-\underline V$.

To control this gap, we additionally impose
\begin{equation*}
\mathcal{R}_y(x^0,y^0)\leq\frac{\mu^4}{L^3},
\end{equation*}
which ensures $c\mathcal{R}_y(x^0,y^0)=O(1)$. Moreover, by the smoothness of $g$ and
$\|u^*(x^0,y^0)\|\leq L_{f,0}/\mu$,
\begin{equation*}
|F(x^0,y^0;u^*(x^0,y^0))-f(x^0,y^0)|
\leq\frac{L_{f,0}}{\mu}
\sqrt{2L_{g,1}\mathcal{R}_y(x^0,y^0)}
=O(\kappa^{-1}).
\end{equation*}
Together with a uniformly bounded initial objective gap
and the lower bound above, this yields
$V^0-\underline V=O(1)$.

This initialization condition is not imposed in the original GALET analysis. Without it, an initialization satisfying $\mathcal{R}_y(x^0,y^0)=\Theta(1)$ contributes $\Theta(c)=\Theta(\kappa^3/\mu)$ to $V^0$, potentially increasing the resulting complexity bound. As discussed in Remark~\ref{initialization_discussion}, this condition can be satisfied by a one-time LL initialization whose cost is independent of the target accuracy. 

Substituting $\alpha=\Theta(\kappa^{-6})$ and $V^0-\underline V=O(1)$ into the descent bound, and using the residual estimates in Appendix~F of \citet{xiao2023generalized}, yields $O(\kappa^6/\epsilon)$ gradient/JVP computations and $O(\kappa^8\epsilon^{-1}\log(\kappa/\epsilon))$ Hessian-vector products.

\vskip 0.2in

\bibliographystyle{plainnat}
\bibliography{nbo}
\end{document}

%% file: macro.tex
\newcommand{\range}[1]{\operatorname{Range}\ (#1)}

\newcommand{\rank}[1]{\operatorname{Rank}\ (#1)}

\newcommand{\dist}[1]{\operatorname{dist}\ (#1)}

\newcommand{\Hk}{\nabla_{22}^2 g(x^k,y^k)}
\newcommand{\He}{\nabla_{22}^2 g(x,y)}
\newcommand{\Hxy}{\nabla_{22}^2 g(x,y)}